\def\a             {\alpha}
\def\Ad            {{\mathrm{Ad}}}

\def\be            {\begin{equation}}
\def\bbC           {\mathbb{C}}

\def\bbR           {\mathbb{R}}
\def\bbT           {\mathbb{T}}
\def\bbZ           {\mathbb{Z}}
\def\bfe           {\mathbf{1}}
\def\can           {\gamma}
\def\canr          {\theta}

\def\cA            {{\mathcal{A}}}

\def\cD            {{\mathcal{D}}}
\def\cE            {{\mathcal{E}}}

\def\cH            {{\mathcal{H}}}

\def\cS            {{\mathcal{S}}}
\def\cT            {{\mathcal{T}}}

\def\cZ            {{\mathcal{Z}}}
\newcommand\co[1]  {\overline{{#1}}}

\def\dim           {{\mathrm{dim}}}

\def\E             {{\mathrm{e}}}
\def\ee            {\end{equation}}
\def\End           {{\mathrm{End}}}
\def\eps           {\varepsilon}

\newcommand\erf[1] {Eq.\ \nolinebreak (\ref{#1})}

\def\ext           {{\mathrm{ext}}}
\newcommand\fig[1] {Fig.\ \nolinebreak \ref{#1}}
\def\furu          {{\mathrm{Furu}}}
\def\Gtwo          {{\mathrm{G}}_2}
\def\Hom           {{\mathrm{Hom}}}
\def\I             {{\mathrm{i}}}
\def\id            {{\mathrm{id}}}

\def\la            {\lambda}
\def\lan           {\langle}
\def\LG            {{\mathit{LG}}}

\def\LIG           {{\mathit{L}}_I{\mathit{G}}}

\def\LSUn          {{\mathit{LSU}}(n)}
\def\LSUz          {{\mathit{LSU}}(2)}

\def\Mat           {{\mathrm{Mat}}}

\def\mult          {{\mathrm{mult}}}
\def\MXN           {{}_M {\mathcal{X}}_N}
\def\MXM           {{}_M {\mathcal{X}}_M}
\def\MXMa          {{}_M^{} {\mathcal{X}}_M^\a}
\def\MXMo          {{}_M^{} {\mathcal{X}}_M^0}

\def\MXMp          {{}_M^{} {\mathcal{X}}_M^+}
\def\MXMm          {{}_M^{} {\mathcal{X}}_M^-}
\def\MXMpm         {{}_M^{} {\mathcal{X}}_M^\pm}
\def\MXMmp         {{}_M^{} {\mathcal{X}}_M^\mp}

\def\MXMppm        {{}_{M_+}^{} {\mathcal{X}}_{M_+}^\pm}

\def\NXN           {{}_N {\mathcal{X}}_N}
\def\NXNd          {{}_N^{} {\mathcal{X}}_N^{\mathrm{deg}}}
\def\NXM           {{}_N {\mathcal{X}}_M}
\def\om            {\omega}
\def\Om            {\Omega}
\def\op            {{\mathrm{opp}}}

\def\ran           {\rangle}

\def\rmD           {{\mathrm{D}}}
\def\rmE           {{\mathrm{E}}}

\def\rmv           {{\mathrm{v}}}
\def\rms           {{\mathrm{s}}}
\def\rmc           {{\mathrm{c}}}
\def\sig           {\sigma}

\def\SLZ           {{\mathit{SL}}(2;\bbZ)}
\def\SOf           {{\mathit{SO}}(5)}

\def\SUd           {{\mathit{SU}}(3)}

\def\SUn           {{\mathit{SU}}(n)}

\def\SUz           {{\mathit{SU}}(2)}

\def\SUf           {{\mathit{SU}}(4)}

\def\tn            {{\tilde{n}}}
\def\tr            {{\mathrm{tr}}}

\def\typei         {type \nolinebreak I}
\def\typeii        {type \nolinebreak II}
\def\typeiii       {type \nolinebreak III}

\def\thinlines{\allinethickness{0.3pt}}
\def\thicklines{\allinethickness{1.0pt}}

\documentclass{fic-l}
\usepackage{amssymb,amsfonts,latexsym,epic,eepic}

\newtheorem{theorem}{Theorem}[section]

\theoremstyle{definition}

\newtheorem{problem}[theorem]{Problem}

\theoremstyle{remark}

\numberwithin{equation}{section}

\newcommand{\abs}[1]{\lvert#1\rvert}

\begin{document}

\title{Modular Invariants and Subfactors}

\author{Jens B\"ockenhauer}
\address{School of Mathematics\\University of Wales Cardiff\\
PO Box 926, Senghennydd Road\\Cardiff CF24 4YH, Wales, UK}
\thanks{Talk given by the first author at the conference
on ``Mathematical Physics in Mathematics and Physics'',
20-25 June 2000, Siena, Italy.}

\author{David E.\ Evans}
\address{School of Mathematics\\University of Wales Cardiff\\
PO Box 926, Senghennydd Road\\Cardiff CF24 4YH, Wales, UK}
\thanks{This project was supported by the EU TMR network
in Non-Commutative Geometry.}

\date{\today}

\dedicatory{This contribution is dedicated to
Sergio Doplicher and John E.\ Roberts\\
on the occasion of their 60th birthdays.}

\begin{abstract}
In this lecture we explain the intimate relationship
between modular invariants in conformal field theory
and braided subfactors in operator algebras.
Our analysis is based on an approach to modular invariants
using braided sector induction (``$\a$-induction'')
arising from the treatment of conformal field theory in the
Doplicher-Haag-Roberts framework.
Many properties of modular invariants which have so far
been noticed empirically and considered mysterious
can be rigorously derived in a very general setting
in the subfactor context.
For example, the connection between modular invariants
and graphs (cf.\ the A-D-E classification for $\SUz_k$)
finds a natural explanation and interpretation.
We try to give an overview on the current state of
affairs concerning the expected equivalence between the
classifications of braided subfactors and modular invariant
two-dimensional conformal field theories.
\end{abstract}

\maketitle

\section{Modular invariants in rational conformal field theory}

It is common knowledge that many, possibly all, models
in rational conformal field theory (RCFT) are related to
current algebras (or ``WZW models'').
A crucial role in the analysis of such current algebra
models is played by their representation theory.
In more mathematical terms: One has to study the
unitary integrable highest weight modules over affine
Lie algebras (cf.\ \cite{K,F}), or,
if you prefer the ``exponentiated version'',
the positive energy representations of loop groups
(cf.\ \cite{PS}), as e.g.\ $\LSUn$.
The positive energy representation are labelled by a level $k$,
a positive integer, and by weights $\la$ in the corresponding
Weyl alcove. Among these, there is a distinguished representation,
the ``vacuum representation'' associated to the weight $\la=0$.
For each positive energy representation
$\pi_\la$, acting on a Hilbert space $H_\la$, we can define
a (specialized) character
\[
\chi_\la(\tau) = \tr_{H_\la} \exp (2\pi\I \tau L_0),
\] 
where $\tau\in\bbC$ is in the upper half plane, and $L_0$
denotes the conformal energy operator which generates the
rotations on the unit circle $S^1$ and which arises from
the affine Lie algebra by the Sugawara construction.
Its lowest eigenvalues $h_\la\ge 0$ is called
``conformal dimension'', and the vacuum representation has
the unique conformal dimension $h_0=0$.
(More generally, the un-specialized characters are defined by
using in addition other variables corresponding to Cartan subalgebra
generators which we omit here for the sake of simplicity.)
It is an important and fascinating fact that the characters
are modular functions. More precisely, at each fixed level $k$
there are unitary matrices $S$ and $T$, the Kac-Peterson matrices
(see \cite{K}), such that
\[
\chi_\la(-1/\tau)= \sum\nolimits_\mu S_{\la,\mu} \chi_\mu(\tau),
\qquad \chi_\la(\tau+1) = \sum\nolimits_\mu
T_{\la,\mu} \chi_\mu(\tau).
\]
Thus there is an action of the modular group $\SLZ$ on
the upper half plane variable $\tau$ with generators
$\cS=\left({0\atop 1}{-1\atop 0}\right)$:
$\tau\mapsto -1/\tau$, and
$\cT=\left({1\atop 0}{1\atop 1}\right)$:
$\tau\mapsto \tau+1$, transforming the family
of characters $\{\chi_\la\}$ unitarily.
Besides the modular relation
\begin{equation}
\label{STcube}
(S T)^3=S^2
\end{equation}
these matrices have
some remarkable properties. In particular, $T$ is diagonal,
and the diagonal entries are up to an overall normalizing
factor given by phases $\exp(2\pi\I h_\la)$,
S is symmetric, $S^2\equiv C$ is a permutation
matrix expressing ``charge conjugation'' which leaves
conformal dimensions invariant since $CT=TC$, and
$S_{\la,0}\ge S_{0,0}>0$. Moreover, $S$ produces
non-negative integers $N_{\la,\mu}^\nu$ by the
Verlinde formula \cite{Ve},
\begin{equation}
\label{verlinde}
N_{\la,\mu}^\nu = \sum_\rho \frac{S_{\la,\rho}}{S_{0,\rho}}
S_{\mu,\rho} S_{\nu,\rho}^*,
\end{equation}
and these integers, called ``fusion rules``, define a
commutative ``fusion rule algebra''.
Combining \erf{STcube} and \erf{verlinde} yields
\[
S_{\la,\mu} = S_{0,0} \sum_{\rho}
\exp(2\pi\I(h_\la+h_\mu-h_\rho)) N_{\la,\mu}^\rho d_\rho,
\]
where $S_{0,0}=1/\sqrt{w}$ with $w=\sum_\la d_\la^2$
and ``quantum dimensions'' $d_\la=S_{\la,0}/S_{0,0}$.

In order to tackle the classification problem of RCFT,
one in particular tries to find all two-dimensional
conformal field theories which contain two copies
of the chiral current algebra as left and right movers
on the (compactified) light rays, and such that
the Hilbert space of the the 2D theory decomposes
upon restriction to the tensor product of chiral
algebras into a direct sum of tensor products
of positive energy representations according to
\begin{equation}
\label{statedecom}
\cH_{\mathrm{phys}} = \bigoplus\nolimits_{\la,\mu}
Z_{\la,\mu} \, H_\la \otimes H_\mu,
\end{equation}
with non-negative integer multiplicities
\begin{equation}
Z_{\la,\mu}=0,1,2,\ldots, \qquad \mbox{and} \qquad Z_{0,0}=1.
\label{constr}
\end{equation}
The latter normalization condition expresses the uniqueness of
the vacuum.

In fact, though there is no direct equivalence, there is a
deep connection between this classification problem and the
problem of classifying modular invariant partition functions
\begin{equation}
\label{Zchar}
\cZ(\tau) = \sum\nolimits_{\la,\mu} Z_{\la,\mu}
\chi_\la(\tau) \chi_\mu(\tau)^*
\end{equation}
such that $\cZ(-1/\tau)=\cZ(\tau)=\cZ(\tau+1)$.
Modular invariant partition functions arise as
continuum limits in statistical mechanics and play
a fundamental role in conformal field theory.
The classification problem for modular invariant
partition functions seems to be more handy than
the classification of 2D RCFT's as it will only
require solutions to the matrix equations $SZ=ZS$,
$TZ=ZT$, subject to the constraints of \erf{constr}.
In fact, as noticed by Gannon \cite{G1}, for given matrices
$S$ and $T$ as above, there are only finitely many
solutions since $w=S_{0,0}^{-2}$ is an overall bound for the
sum of all entries of $Z$.
There is also an inequality
\begin{equation}
\label{Zbound}
Z_{\la,\mu} \le d_\la d_\mu
\end{equation}
for each individual entry of a modular invariant
coupling matrix \cite{BE5}. So clearly for a fixed model,
e.g.\ for a certain $\SUn$ at a fixed level $k$, there is
only a finite number of modular invariants. However,
ambitious scientists are usually interested in classifying
modular invariants for entire series of models, e.g.\
$\SUn$, fixed rank but all levels.
Unfortunately, this ambition turned out to meet hard problems,
and complete classifications exist only for $n=2$,
the celebrated A-D-E classification of \cite{CIZ1,CIZ2,Kt},
and for $n=3$ \cite{G2}.
Gannon has recently informed us that he has completed
the $\SUf$ case up to levels $k\le 5,000$ and with
similar bounds also other Lie groups, but otherwise
there are still many open problems around.
(At fixed low levels $k<4$, there are however
classifications for all $\SUn$, see \cite{G3}.)
Anyway, the classification problem of modular invariants
is definitely a fascinating area of mathematical physics
with many deep connections to other branches of mathematics
and physics, see e.g.\ \cite{G5}.

To be a bit more illustrative, let us see some examples.
For $SU(2)$ at level $k$, the positive energy representations
are labelled simply by spins $\la=0,1,2,....,k$.
At level $k=6$, there are two modular invariants,
which read, when written in the form of \erf{Zchar},
\[
\begin{split}
\cZ_{\mathrm{A}_7} &= \abs{\chi_0}^2 + \abs{\chi_1}^2
+ \abs{\chi_2}^2 + \abs{\chi_3}^2 + \abs{\chi_4}^2
+ \abs{\chi_5}^2 + \abs{\chi_6}^2, \\
\cZ_{\mathrm{D}_5} &= \abs{\chi_0}^2 + \abs{\chi_2}^2
+ \abs{\chi_4}^2 + \abs{\chi_6}^2
+ \chi_1\chi_5^* + \chi_3\chi_3^* + \chi_5\chi_1^*.
\end{split}
\]
At level $k=16$ there are three:
\[
\begin{split}
\cZ_{\mathrm{A}_{17}} &= \sum\nolimits_{\la=0}^{16} \abs{\chi_\la}^2 \\
\cZ_{\mathrm{D}_{10}} &=  \abs{\chi_0 + \chi_{16}}^2
+ \abs{\chi_2 + \chi_{14}}^2 + \abs{\chi_4 + \chi_{12}}^2
+ \abs{\chi_6 + \chi_{10}}^2 + 2\abs{\chi_8}^2, \\
\cZ_{\mathrm{E}_7} &= \abs{\chi_0 + \chi_{16}}^2
+ \abs{\chi_4 + \chi_{12}}^2 + \abs{\chi_6 + \chi_{10}}^2
+ \abs{\chi_8}^2 \\
&\qquad\qquad + (\chi_2 + \chi_{14})\chi_8^*
+ \chi_8 (\chi_2 + \chi_{14})^*.
\end{split}
\]
There is so much structure visible already in these relatively
simple examples that we would like to comment on this.
Let us first explain the labelling of the $\cZ$'s by A-D-E Dynkin
diagrams. It was noticed in \cite{CIZ1} that the diagonal
terms of each invariant appearing at level $k$ reproduce
exactly the Coxeter exponents of one of the Dynkin diagrams
with Coxeter number $k+2$, and this amounts to a bijective
correspondence between all the $\SUz_k$, $k=1,2,...$,
modular invariants
(The list of \cite{CIZ1} was proven to be complete in \cite{CIZ2,Kt},
and for a more recent elegant proof see \cite{G4})
and all the A-D-E Dynkin diagrams.
More precisely, the labelling the invariants by Dynkin diagrams
is such that the diagonal entry $Z_{\rho,\rho}$ of the invariant
associated to the A-D-E diagram $G_1$ is exactly the multiplicity
of the eigenvalue
\[
\frac{S_{1,\rho}}{S_{0,\rho}} = 2\cos \left(
\frac{(\rho+1)\pi}{k+2} \right)
\]
of the (adjacency matrix of the) graph $G_1$.
This empirical observation was certainly awaiting
a good explanation! Nahm found \cite{N} a systematic
relation between the diagonal part of the $\SUz$ invariants
and Lie algebra exponents using quaternionic coset spaces.
But Nahm's ideas do not explain that the A-D-E observation
turned out to be just the ``$\SUz$-spin-1-tip of the iceberg'':
First of all, the A-D-E diagrams $G_1$ are just the
spin-1 member of a family of graphs $G_\la$, $\la=0,1,2,...,k$,
which form a non-negative integer valued matrix representation
(nimrep, for short) of the $\SUz_k$ Verlinde fusion rules rules:
\[
G_\la G_\mu = \sum\nolimits_\nu N_{\la,\mu}^\nu G_\nu ,
\]
with eigenvalue multiplicities
$\mult_{G_\la}(S_{\la,\rho}/S_{0,\rho})=Z_{\rho,\rho}$.
(For the Dynkin diagrams of type A, the $G_\la$'s are
just the fusion matrices, and the appearance of each
character $\gamma_\rho(\cdot)=S_{\cdot,\rho}/S_{0,\rho}$
with multiplicity $Z_{\rho,\rho}=1$ for all $\rho$ is
just the Verlinde formula.)
And even more, guided by the observations for $\SUz$,
Di Francesco and Zuber found \cite{DZ1,DZ2,DiF}
(see also related work \cite{PZ,BPPZ}) that there are
graphs and nimreps for $\SUd$ and also higher rank $\SUn$
modular invariants which fall into the analogous scheme,
i.e.\ such that you just have to replace the labels (``weights''),
fusion rules and S-matrices by the corresponding $\SUn$ data.
In the subfactor context, we will see that these
nimreps arise from a certain braided sector induction,
called ``$\a$-induction''; the matrix entries are non-negative
integers since they are dimensions of certain intertwiner spaces,
the representation property is due to the fact that $\a$-induction
preserves fusion rules, and a general theorem determining the
character multiplicities to be given by the diagonal entries
of the modular invariant coupling matrix can be proven.

Next we address the distinction of \typei\ and \typeii\ modular
invariants. Note that, whatever model we are looking at,
there will always be at least one solution,
the diagonal partition function
\[
\cZ= \sum\nolimits_\la |\chi_\la|^2,
\]
which is always modular invariant, equivalently expressed in
the fact that the unit matrix, $Z_{\la,\mu}=\delta_{\la,\mu}$,
always commutes with $S$ and $T$. For $\SUz$, these are
the invariants labelled by Dynkin diagrams A$_{k+1}$.
More generally, there may be permutation invariants
\[
\cZ = \sum\nolimits_\la \chi_\la \chi_{\vartheta(\la)}^* \,,
\]
whenever $\vartheta$ is a permutation of the labels which preserves
the fusion rules, the vacuum, and the conformal dimensions
modular integers. The above displayed D$_5$ invariant for
$\SUz_6$ is an example for such an automorphism.
In general the charge conjugation matrix $C$ is also such
an automorphism --- which is however trivial in the special
case of $\SUz$.
Moore and Seiberg argue in \cite{MS2} (see also \cite{DV})
that after a ``maximal extension of the chiral algebra''
(the hardest part is to make this mathematically precise)
the partition function of a RCFT is at most a permutation matrix
$Z^\ext_{\tau,\tau'}=\delta_{\tau,\vartheta(\tau')}$, where
$\tau,\tau'$ now label the representations of the extended
chiral algebra and $\vartheta$ denotes a permutation of these
with analogous invariance properties.
Decomposing the extended characters $\chi^\ext_\tau$
in terms of the original characters $\chi_\la$, we have
$\chi^\ext_\tau=\sum_\la b_{\tau,\la} \chi_\la$ for some
non-negative integral branching coefficients $b_{\tau,\la}$.
The maximal extension yields the coupling matrix expression
\begin{equation}
\label{extaut}
Z_{\la,\mu}=\sum\nolimits_\tau
b_{\tau,\la} b_{\vartheta(\tau),\mu}.
\end{equation}
Di Francesco and Zuber \cite{DZ2} called invariants which
arise from the diagonal invariant of the maximal extension
(i.e.\ for which $\vartheta$ is trivial) ``\typei'',
and invariants corresponding to non-trivial automorphisms
of the extended fusion rules were called ``\typeii''.
Looking at the above displayed $\SUz_6$ and $\SUz_{16}$
invariants, we find for example that
$\cZ_{\mathrm{A}_7}$, $\cZ_{\mathrm{A}_{17}}$ and
$\cZ_{\mathrm{D}_{10}}$ are \typei\ whereas
$\cZ_{\mathrm{D}_5}$ and $\cZ_{\mathrm{E}_7}$
are \typeii.

Let us finally remark that the A-D-E classification generalizes
in some sense to the entire classification problem of modular
invariants in RCFT. The class of diagonal modular invariants and
their conjugations is often denoted by $\cA$.
A wider class is given by the ``simple current invariants''
(see e.g.\ \cite{SY1,SY2,GRS}) for which $Z_{\la,\mu}\neq0$
implies $N_{\sigma,\la}^\mu=1$ for some simple current $\sigma$
(i.e.\ a label with $d_\sigma=1$). The class of simple current
invariants minus the $\cA$ class is often denoted by $\cD$.
The remaining modular invariants, which are typically relatively
few, are called ``exceptionals'' and their class is denoted
by $\cE$. In fact, considering a loop group model
(with fixed rank) at all levels, the $\cA$ and, if there
are simple currents, $\cD$ classes give infinite, very
well-behaving series of invariants and there is only a finite
number of exceptionals. The graphs which have been associated
to diagonal $\cA$ invariants are basically the Weyl alcove
with edges corresponding to the fusion with the fundamental
generator. The graphs which have been associated to
simple current $\cD$ invariants are simply orbifolds of
the $\cA$-type graphs with respect to a cyclic simple
current symmetry. (E.g.\ the Dynkin diagrams D$_{\varrho+2}$
are $\bbZ_2$ orbifolds of the A$_{2\varrho+1}$ graphs at
even levels $k=2\varrho$, $\varrho=2,3,...$.)
The graphs associated to $\cE$ invariants
may be considered as orbifolds with
respect to a more subtle, non-group-like symmetry.
In any case, it was noticed (see \cite{DiF}) that
one can associate intrinsic fusion rules algebras
to the graphs exactly in the \typei\ cases, and that
they contain fusion subalgebras corresponding to the
Verlinde fusion rules of the extended characters.
All this finds a natural explanation in the subfactor
framework.

\section{Modular invariants from subfactors through $\a$-induction}

Now let us switch to a different topic: Subfactors.
(See \cite{EK} as a general reference on subfactors.)
At the first sight, this topic does not seem to have
anything to do with modular invariants.
However, soon we shall see that it does,
and in fact quite a lot!

Let $A$ and $B$ be \typeiii\ von Neumann factors.
A unital $\ast$-homomorphism $\rho:A\rightarrow B$
is called a $B$-$A$ morphism. The positive number
$d_\rho=[B:\rho(A)]^{1/2}$ is called the
statistical dimension of $\rho$; here $[B:\rho(A)]$ is
the Jones-Kosaki index \cite{J1,Ko} of the subfactor
$\rho(A)\subset B$. Some $B$-$A$ morphism $\rho'$
is called equivalent to $\rho$ if $\rho'=\Ad(u)\circ\rho$
for some unitary $u\in B$. The equivalence class $[\rho]$
of $\rho$ is called the $B$-$A$ sector of $\rho$.
If $\rho$ and $\sig$ are $B$-$A$
morphisms with finite statistical dimensions, then
the vector space of intertwiners
\[ \Hom(\rho,\sig)=\{ t\in B: t\rho(a)=\sig(a)t \,,
\,\, a\in A \}  \]
is finite-dimensional, and we denote its dimension by
$\lan\rho,\sig\ran$.
In fact $\lan\rho,\sig\ran\le d_\rho d_\sigma$.
A $B$-$A$ morphism is called irreducible if
$\langle \rho,\rho \rangle=1$, i.e.\ if
$\Hom(\rho,\rho)=\bbC\bfe_B$.
Then, if $\langle\rho,\tau\rangle\neq 0$ for some
(possibly reducible) $B$-$A$ morphism $\tau$,
then $[\rho]$ is called an irreducible subsector of
$[\tau]$ with multiplicity $\langle\rho,\tau\rangle$.
An irreducible $A$-$B$ morphism $\co\rho$ is a conjugate
morphism of the irreducible $\rho$ if and only if
$[\co\rho\rho]$ contains the trivial sector $[\id_A]$
as a subsector, and then
$\langle\rho\co\rho,\id_B\rangle=1=\langle\co\rho\rho,\id_A\rangle$
automatically \cite{I1}.
The map $\phi_\rho:B\rightarrow A$,
$b\mapsto r_\rho^* \co\rho(b)r_\rho$, is called the
(unique) standard left inverse and satisfies
\[
\phi_\rho(\rho(a)b\rho(a'))=a\phi_\rho(b)a' \,,
\quad a,a'\in A\,,\quad b\in B \,.
\]

Now let $N$ be a \typeiii\ factors.
We assume that we have a given finite system of $N$-$N$ morphisms
$\NXN$, i.e.\ $\NXN\subset\End(N)$ is finite and
\begin{itemize}
\item each $\la\in\NXN$ is irreducible,
  i.e.\ $\la(N)'\cap N=\bbC\bfe_N$,
\item each $\la\in\NXN$ has finite statistical dimension,
  i.e.\ $d_\la=[N:\la(N)]^{1/2}<\infty$,
\item the morphisms are pairwise inequivalent,
  i.e.\ $\langle\la,\mu\rangle=0$ whenever $\la\neq\mu$,
\item the identity morphism belongs to the system,
  i.e.\ $\id\in\NXN$,
\item the system is closed under conjugation,
  i.e.\ for each $\la\in\NXN$ there is some $\co\la\in\NXN$
  such that $\langle\co\la\la,\id_N\rangle=1$,
\item the system is closed under fusion,
  i.e.\ any composition $[\la]\cdot[\mu]\equiv[\la\circ\mu]$,
  $\la,\mu\in\NXN$ has only subsectors arising from $\NXN$;
  we write
  $[\la]\cdot[\mu]=\bigoplus_{\nu\in\NXN}
  \langle\la\mu,\nu\rangle [\nu]$.
\end{itemize}

We now also assume that our system $\NXN$ is braided:
For any pair $\la,\mu\in\NXN$ there is a unitary operator
$\eps^+(\la,\mu)\in\Hom(\la\mu,\mu\la)$ such that the
braiding fusion equations,
\begin{equation}
\begin{split}
\rho(t) \eps^+(\la,\rho) 
&=\eps^+(\mu,\rho) \mu(\eps^+(\nu,\rho))t ,\\
t \eps^+(\rho,\la)
&=\mu (\eps^+(\rho,\nu)) \eps^+(\rho,\mu)\rho (t) , \\
\rho(t)^* \eps^+(\mu,\rho) \mu( \eps^+(\nu,\rho))
&= \eps^+(\la,\rho) t^* , \\
t^* \mu( \eps^+(\rho,\nu)) \eps^+(\rho,\mu)
&=  \eps^+ (\rho,\la) \rho(t)^* ,
\end{split}
\label{BFE}
\end{equation}
hold whenever $\la,\mu,\nu\in\NXN$ and
$t\in \Hom (\lambda,\mu\nu)$.
The unitaries $\eps^+(\la,\mu)$ are called statistics operators.
Note that a braiding $\eps\equiv\eps^+$ always
comes along with another ``opposite'' braiding $\eps^-$,
namely operators $\eps^-(\la,\mu) = \eps^+(\mu,\la)^*$,
satisfy the same relations. The unitaries $\eps^+(\la,\mu)$
and $\eps^-(\la,\mu)$ are different in general but may coincide
for some $\la$, $\mu$ --- at least if one of them is the identity
morphism since \erf{BFE} implies $\eps^+(\la,\id)=\bfe$. 

Exactly as in the DHR theory (cf.\ \cite{H}),
$d_\la \phi_\la(\eps^+(\la,\la))\in\Hom(\la,\la)$
is unitary, and so one defines the statistics phase
$\om_\la\in\bbT$ by
\[
d_\la \phi_\la(\eps^+(\la,\la))\in\Hom(\la,\la)
 = \om_\la \bfe.
\]
Now consider the following matrices $\Om$ and $Y$,
\[
\begin{split}
\Om_{\la,\mu} &= \delta_{\la,\mu} \om_\la ,\\
Y_{\la,\mu} &= \sum\nolimits_{\rho\in\NXN}
\frac{\om_\la \om_\mu}{\om_\rho} \langle\la\mu,\rho\rangle
d_\rho ,
\end{split}
\]
with indices labelled by $\NXN$, and we will use the
label ``0'' for the identity morphism $\id\in\NXN$.
Then one checks that $Y$ is symmetric,
that $Y_{\co\la,\mu}=Y_{\la,\mu}^*$ as well as
$Y_{\la,0}=d_\la$.
The Y- and $\Omega$-matrices obey
$\Omega Y \Omega Y \Omega = z Y$ where
$z=\sum_\la d_\la^2 \om_\la$ \cite{R1,FG1,FRS2}.
If $z\neq 0$ we put $c=4\arg(z)/\pi$, which is defined
modulo 8, and call it the ``central charge'', and then
statistics S- and T-matrices defined by
\[ S = |z|^{-1} Y \,, \qquad T = \E^{-\I\pi c/12} \Omega \]
hence fulfill $TSTST=S$.
With $C$ denoting the sector conjugation matrix,
$C_{\la,\mu}=\delta_{\la,\co\mu}$, one finds
$CT=TC$ and $CS=SC$.
Rehren showed in \cite{R1} that the following
conditions are equivalent:
\begin{itemize}
\item The braiding $\eps$ is non-degenerate, i.e.\
  $\eps^+(\la,\mu)=\eps^-(\la,\mu)$ for all $\mu\in\NXN$
  only if $\la=\id$.
\item One has $|z|^2=w$ with the global index
$w=\sum_\la d_\la^2$ and $S$ is unitary, so that
$S$ and $T$ are indeed the standard generators
$\cS=\left({0\atop 1}{-1\atop 0}\right)$ and
$\cT=\left({1\atop 0}{1\atop 1}\right)$
in a unitary representation of $\SLZ$,
i.e.\ fulfill \erf{STcube}, and indeed $S^2=C$.
Moreover, a Verlinde formula holds (cf.\ \erf{verlinde}):
\[
\sum_{\rho\in\NXN} \frac{S_{\la,\rho}}{S_{0,\rho}}
S_{\mu,\rho} S_{\nu,\rho}^* = \langle\la\mu,\nu\rangle.
\]
\end{itemize}
So here we have obtained another representation of
the modular group $\SLZ$, and in fact we seem to be
dealing with precisely the same categorical structures as
in RCFT (at least if the braiding is non-degenerate)!

Now let $\LG$ be a loop group (associated to a simple,
simply connected loop group $G$). Let $\LIG$ denote
the subgroup of loops which are trivial off some
proper interval $I\subset S^1$. Then in each level $k$
vacuum representation $\pi_0$ of $\LG$, we naturally obtain
a net\footnote{In fact a proper net is obtained only if
we remove a ``point at infinity'' from the circle $S^1$.}
of \typeiii\ factors $\{N(I)\}$ indexed by
proper intervals $I\subset S^1$ by taking
$N(I)=\pi_0(\LIG)''$ (see \cite{W,FG2,Bg}).
Since the DHR selection criterion (cf.\ \cite{H}) is met in
the (level $k$) positive energy representations $\pi_\la$,
there are DHR endomorphisms $\la$ naturally associated with them.
(By some abuse of notation we use the same symbols for
labels and endomorphisms.)
Now it is very natural to expect the following
(a conjecture which actually goes beyond loop groups,
see e.g.\ \cite{FG2}): We anticipate that the statistics
phases are the exponentiated conformal dimensions,
i.e.\ that
\begin{equation}
\om_\la=\exp(2\pi\I h_\la),
\label{spst}
\end{equation}
and that the RCFT Verlinde fusion coincides with
the (DHR superselection) sector fusion, i.e.\ that
\begin{equation}
N_{\la,\mu}^\nu=\langle\la\mu,\nu\rangle .
\label{RCFT=DHR}
\end{equation}
(And in turn that the RCFT quantum dimensions
equal the statistical dimensions.)
In other words, we expect that the normalized
matrices $Y$ and $\Om$, namely the statistics
S- and T-matrices are identical with the Kac-Peterson
S- and T-modular matrices which perform the
conformal character transformations.
Fortunately general results have been proven \cite{FG1,FRS2,GL2}
for the ``conformal spins statistics theorem'', \erf{spst}.
For \erf{RCFT=DHR} there are proofs available
\cite{W,Loke,Tol,B1,B2} unfortunately only for special
models.\footnote{Antony Wassermann has informed us that
he has extended his results for $\SUn_k$ fusion \cite{W}
to all simple, simply connected loop groups; and with
Toledano-Laredo all but E$_8$ using a variant of the
Dotsenko-Fateev differential equation considered in
his thesis \cite{Tol}.}
Anyway, we will be mainly concerned with $\SUn_k$ here,
so we can take the equality of statistics and Kac-Peterson
matrices for granted, thanks to \cite{W}.

To summarize the above paragraphs, we have seen that
a factor with a (non-degenerately braided system of
endomorphisms gives rise to a unitary representation
of the modular group $\SLZ$ via matrices $S$ and $T$
which are analogous to the Kac-Peterson matrices in RCFT.
So what about modular invariants? As we shall see,
modular invariants appear naturally in the operator
algebraic setting when we consider \textit{sub}factors
with a braiding. Suppose we have an embedding of our
factor $N$ in a larger factor $M$, i.e.\ we have a
subfactor $N\subset M$. Let
$\iota:N\hookrightarrow M$ be the inclusion map which we
may consider as an $M$-$N$ morphism. Choose a representative
$\co\iota:M\rightarrow N$ of the conjugate $N$-$M$ sector.
Then $\canr=\co\iota\iota$ is Longo's dual canonical
endomorphism, and we call $N\subset M$ a braided subfactor
if its sector $[\canr]$ decomposes exclusively into
subsectors of our braided system $\NXN$, i.e.\ if
\[
[\canr] =\bigoplus_{\rho\in\NXN} n_\rho [\rho],
\qquad n_\rho=\langle\rho,\canr\rangle \in \{ 0,1,2,\ldots \}.
\]
It is straightforward to extend a braiding $\eps$ on $\NXN$
to the set $\Sigma(\NXN)$ of all to equivalent morphisms and
direct sums (see e.g.\ \cite{BEK1}). Then one can define
the $\a$-induced morphisms $\a^\pm_\la\in\End(M)$
for $\la\in\Sigma(\NXN$) by the Longo-Rehren formula \cite{LR}
which concretely realizes a ``cohomological extension''
suggested by Roberts \cite{Ro} about 24 years ago.
Namely one puts
\[
\a_\la^\pm = \co\iota^{\,-1} \circ \Ad
(\eps^\pm(\lambda,\canr)) \circ \lambda \circ \co\iota .
\]
Then $\a^+_\la$ and $\a^-_\la$ extend $\la$, i.e.\
$\a^\pm_\la\circ\iota=\iota\circ\la$, which in turn implies
$d_{\a_\la^\pm}=d_\la$ by the multiplicativity of
the minimal index \cite{L3}. 
Let $\can=\iota\co\iota$ denote Longo's canonical
endomorphism \cite{L1} from $M$ into $N$.
Then there is an isometry $v\in\Hom(\id,\can)$
such that any $m\in M$ is uniquely decomposed
as $m=nv$ with $n\in N$ \cite{L2}.
Thus the action of the extensions
$\a^\pm_\la$ are uniquely characterized by the relation
$\a^\pm_\la(v)=\eps^\pm(\la,\canr)^* v$ which can be
derived from the braiding fusion equations \erf{BFE}.
We have $\a_{\co\la}^\pm$ is a conjugate for $\a_\la^\pm$,
moreover $\a_{\la\mu}^\pm=\a_\la^\pm \a_\mu^\pm$
if also $\mu\in\Sigma(\NXN)$, and clearly
$\a_{{\rm{id}}_N}^\pm={{\rm{id}}}_M$
(proofs can be found in \cite{BE1}, and for a similar
framework -- the relations are explained in \cite{X2} -- in
the earlier work \cite{X1}). In general one has
\[ \Hom(\la,\mu) \subset \Hom(\a^\pm_\la,\a^\pm_\mu)
\subset \Hom(\iota\la,\iota\mu) \,,\qquad\la,\mu\in\Sigma(\NXN). \]
Now let us count the common subsectors of $\a$-induced
morphisms with different chirality ``$+$'' and ``$-$''
by defining a ``coupling matrix'' $Z$ with entries
\begin{equation}
\label{Zdef}
Z_{\la,\mu} = \langle \a^+_\la,\a^-_\mu \rangle,
\qquad \la,\mu\in\NXN .
\end{equation}
Clearly, $Z_{\la,\mu}$ are non-negative integers since
dimensions of vector spaces. Moreover, $Z_{0,0}=1$ due
to $\a_{{\rm{id}}_N}^\pm={{\rm{id}}}_M$.
Also note that \erf{Zdef} immediately yields \erf{Zbound}.
It has been shown in \cite[Thm.\ 5.7]{BEK1} that in fact
\[
YZ=ZY, \qquad \Om Z=Z\Om ,
\]
no matter whether the braiding is degenerate or not.
So here we have a notion of modular invariants
arising from subfactors which extends the non-degenerate
(i.e.\ modular) case to non-unitary S-matrices.
We now would like to see more structure and to
understand e.g.\ the connection between modular
invariants and graphs or the Moore-Seiberg machinery
involving fusion rules automorphism invariants,
\typei\ and \typeii\ invariants etc.\ from the
operator algebraic point of view. As we shall see,
this viewpoint opens up new insights and resolves
so far somewhat mysterious phenomena. For that,
we have to analyze the structure of $\a$-induced sectors.

\section{Structure of modular invariants from subfactors:\\
 Induced sector systems, fusion and graphs}

Let $\MXM\subset\End(M)$ denote a system of endomorphisms
consisting of a choice of representative endomorphisms of
each irreducible subsector of sectors of the form
$[\iota\la\co\iota]$, $\la\in\NXN$.
We choose $\id\in\End(M)$
representing the trivial sector in $\MXM$.
Then we define similarly the ``chiral'' systems
$\MXMpm$ and the ``$\a$-system'' $\MXMa$ to be the
subsystems of endomorphisms $\beta\in\MXM$
such that $[\beta]$ is a subsector of $[\a^\pm_\la]$
and of of $[\a_\la^+\a_\mu^-]$, respectively,
with $\la,\mu\in\NXN$ varying.
(Any subsector of $[\a_\la^+\a_\mu^-]$ is
automatically a subsector of  $[\iota\nu\co\iota]$
for some $\nu\in\NXN$.)
The ``neutral'' (or ``ambichiral'') system is defined
as the intersection $\MXMo=\MXMp\cap\MXMm$, so that
$\MXMo \subset \MXMpm \subset \MXMa \subset \MXM$.
Their ``global indices'', i.e.\ their sums over the
squares of the statistical dimensions are denoted
by $w_0$, $w_\pm$, $w_\a$ and $w$, and thus fulfill
$1 \le w_0 \le w_\pm \le w_\a \le w$.
It turns out that the relative sizes
of these systems, which are measured by such global indices,
are completely encoded in the coupling matrix $Z$,
namely \cite[Prop.\ 3.1]{BEK2}
\begin{equation}
\label{chirglob}
w_+ = \frac w{\sum_{\la\in\NXN} d_\la Z_{\la,0}}
= \frac w{\sum_{\la\in\NXN} Z_{0,\la} d_\la} = w_- ,
\end{equation}
and \cite[Prop.\ 3.1]{BE4}
\begin{equation}
\label{aglob}
w_\a = \frac w{\sum_{\la\in\NXNd} Z_{0,\la} d_\la} ,
\qquad w_0 = \frac{w_+^2}{w_\a} ,
\end{equation}
where $\NXNd\subset\NXN$ denotes the subsystem of
degenerate morphisms. As a corollary of \erf{aglob} one
obtains that non-degeneracy of the braiding
(i.e.\ $\NXNd=\{\id\}$) implies the
``generating property'' $\MXMa=\MXM$.

Though $\NXN$ is braided by assumption, the systems
$\MXMpm$ or $\MXM$ will in general not be.
In fact, as constructed explicitly in \cite{BE3},
there is only a relative braiding between $\MXMp$
and $\MXMm$, and this restricts to a proper braiding
on the intersection $\MXMo$. However, the systems
$\MXMpm$ can even be non-commutative.
A criterion was found for the case
of a non-degenerate braiding: Namely,
if we consider the fusion rules of $\MXMpm$ as
finite-dimensional $C^*$-algebras $\furu(\MXMpm)$,
the we have \cite[Thm.\ 4.11]{BEK2}
\begin{equation}
\label{chirfudec}
\furu(\MXMpm) \simeq \bigoplus_{\la\in\NXN}
\bigoplus_{\tau\in\MXMo} \Mat(b^\pm_{\tau,\la})
\end{equation}
with ``chiral branching coefficients''
$b^\pm_{\tau,\la}=\langle \tau,\a^\pm_\la\rangle$.
The analogous result for $\MXM$ reads \cite[Thm.\ 6.8]{BEK1}
(also provided that the braiding on $\NXN$ is non-degenerate)
\begin{equation}
\label{fullfudec}
\furu(\MXM) \simeq \bigoplus_{\la,\mu\in\NXN}
\Mat(Z_{\la,\mu}).
\end{equation}
(The latter decomposition was claimed in a similar form
in \cite{O} in the context of Goodman-de la Harpe-Jones
subfactors related to $\SUz_k$ modular invariants.)
Equivalently one may determine the irreducible decomposition
of the ``regular representations'' $\pi^\pm_{\mathrm{reg}}$,
$\pi_{\mathrm{reg}}$ of the fusion rule algebras
$\furu(\MXMpm)$, $\furu(\MXM)$, respectively, and the
corresponding irreducible decompositions then read
\begin{equation}
\label{regreps}
\begin{split}
\pi^\pm_{\mathrm{reg}} &\simeq \bigoplus_{\la\in\NXN}
\bigoplus_{\tau\in\MXMo} b^\pm_{\tau,\la} \pi^\pm_{\tau,\la},\\
\pi_{\mathrm{reg}} & \simeq \bigoplus_{\la,\mu\in\NXN}
Z_{\la,\mu} \pi_{\la,\mu},
\end{split}
\end{equation}
with multiplicities given by the dimensions of the
irreducible representations. Note that this is essentially
the block-diagonalization of the fusion matrices of the
intrinsic fusion rules of the systems $\MXMpm$ and $\MXM$.
Reflecting commutativity properties of the induced
sectors $[\a^+_\nu]$ ($[\a^-_\nu]$), these are scalars
in the irreducible representations $\pi^+_{\tau,\la}$
($\pi^-_{\tau,\mu}$) and $\pi_{\rho,\sigma}$, given by
$S_{\nu,\la}/S_{0,\la}$ ($S_{\nu,\mu}/S_{0,\mu}$)
and $S_{\nu,\rho}/S_{0,\rho}$ ($S_{\nu,\sigma}/S_{0,\sigma}$),
respectively \cite[Cor.\ 4.15]{BEK2}.
However, we can now easily see that the entire system $\MXM$
has non-commutative fusion if and only if an entry of $Z$ is
strictly larger than one. There are lots of examples, the
simplest given by the D$_{\mathrm{even}}$ series of $\SUz_k$.
Similarly we see that the chiral system $\MXMpm$ have
non-commutative fusion if and only if there is a
chiral branching coefficient $b^\pm_{\tau,\la}$
strictly larger than one.
To find examples one has to dig a bit further.
In fact this happens for a series of conformal inclusions
$\SUn_n\subset{\mathit{SO}}(n^2-1)_1$ for $n\ge 4$,
and indeed a non-commutative chiral fusion structure was
found for the case $n=4$ by direct computation in \cite{X1}
(see also \cite{BE2} for a treatment using the Longo-Rehren
$\a$-induction).

So how can we interpret non-commutative fusion rules?
Why is there a relative braiding between the possibly
non-commutative chiral systems $\MXMp$ and $\MXMm$?
For this we should think of our endomorphisms in $\NXN$
again as DHR endomorphisms of a whole net of algebras
$\{N(I)\}$ over the punctured circle rather
than of a single local algebra $N(I_\circ)$.
In this context the (subsectors of the) $\a$-induced sectors
are in fact solitonic, i.e.\ left or right half-line localized
depending on their $\pm$-chirality, and their respective
localization regions intersect exactly on the chosen
interval $I_\circ$ (see \cite{LR}).
Then in the DHR framework such solitonic or ``twisted'' sectors
of different chirality can be ``pulled apart'' and commuted,
a procedure which provides in fact unitaries obeying partial
braiding properties. However, for sectors with the same
half-line localization (or even without any localization
as is the case for sectors in the mixed system $\MXMa$)
this procedure does not work, and so there is no reason why
such sectors should have commutative fusion.
The neutral system $\MXMo$ however corresponds to proper DHR
endomorphisms, and their DHR statistics operators are precisely
their restricted relative braiding \cite[Prop.\ 3.15]{BE3}.

But let us return to the regular representations in \erf{regreps}.
The representation theoretic point of view has the advantage
that we can also consider the left multiplication of $M$-$M$
sectors on $M$-$N$ sectors: Let $\MXN$ denote a system
consisting of a choice of representative $M$-$N$ morphisms
of irreducible subsectors of sectors $[\iota\la]$, with
$\la\in\NXN$ varying. Then $\MXN$ has no intrinsic fusion
structure, nevertheless we can consider the representation
$\varrho$ of $\furu(\MXM)$ arising from multiplication of
$\MXM$ on $\MXN$. Its representation matrices $\varrho(\beta)$,
$\beta\in\MXM$ are given by
\[
\varrho(\beta)_{a,b} = \langle b, \beta a \rangle ,
\qquad a,b \in \MXN .
\]
The irreducible decomposition of $\varrho$ has been
determined in \cite[Thm.\ 6.12]{BEK1} to be
\begin{equation}
\label{varrhorep}
\varrho \simeq \bigoplus_{\la\in\NXN} \pi_{\la,\la}.
\end{equation}
Now consider the following matrix representation of
$\furu(\NXN)$, with representation matrices $G_\la$,
$\la\in\NXN$, with non-negative integer entries
\begin{equation}
\label{Gmatrix}
(G_\la)_{a,b} = \langle b, \a^\pm_\la a \rangle ,
\qquad a,b \in \MXN ,
\end{equation}
i.e.\ $G_\la=\sum_{\beta\in\MXMpm}
\langle\beta,\a^\pm_\la\rangle \varrho(\beta)$.
(It does not depend on the choice of $\pm$-chirality.)
Thanks to \erf{varrhorep}, we now know the complete
diagonalization of the $G_\la$'s: The eigenvalues
are of the form $S_{\la,\rho}/S_{0,\rho}$, and
the multiplicities are given by the dimensions
$\dim(\pi_{\rho,\rho})=Z_{\rho,\rho}$, the diagonal
entries of the modular invariant \cite[Thm.\ 4.16]{BEK2}.
So here we have obtained systematically a
nimrep of the original Verlinde fusion algebra
$\furu(\NXN)$ from a subfactor, and it has spectrum
canonically associated with the diagonal part of the
corresponding modular invariant.
Trivially, the non-negative integer valued matrices can
be read as adjacency matrices of graphs and this
way we obtain in particular graphs associated to some
modular invariant.

Let us take a closer look at the $M$-$N$ system $\MXN$.
For each $a\in\MXN$ we can choose a conjugate $N$-$M$
morphism $\co a$, and this way we obtain a system $\NXM$
of irreducibles $N$-$M$ morphisms such that in particular
$\co\iota\in\NXM$. Now for any such $\co a\in\NXM$
there is an irreducible subfactor $\co a(M)\subset N$
and we can form the Jones basic extension
\[
\co a(M) \subset N \subset M_a .
\]
For the injection homomorphism $\iota_a:N\hookrightarrow M_a$,
we can choose a conjugate $M_a$-$N$ morphism ${\co\iota}_a$
such that ${\co\iota}_a(M_a)=\co a(M)$. Hence we have an
isomorphism $\varphi_a:M_a\rightarrow M$ given by
$\varphi_a={\co a}^{\,-1}\circ{\co\iota}_a$ with inverse
$\varphi_a^{-1}={\co\iota}_a^{\,-1}\circ\co a$.
Note that the dual canonical endomorphism
$\canr_a={\co\iota}_a \iota_a$ of the subfactor $N\subset M_a$
can also be written as $\canr_a=\co a a$, and that the
Jones index is $[M_a:N]=d_{\canr_a}=d_a^2$.
Next we  consider $\a$-induction of $\la\in\Sigma(\NXN)$
for $N\subset M_a$:
\[
\a^\pm_{a;\la} = {\co\iota}_a^{\,-1} \circ
\Ad(\eps^\pm(\la,\canr_a) \circ \la \circ {\co\iota}_a .
\]
It follows from \cite[Props.\ 3.1,3.3]{BEK1} that
$\eps^\pm(\la,\canr_a)\equiv\eps^\pm(\la,\co a a)$ can be
written as $\eps^\pm(\la,\canr_a)=\co a(U_\la^\pm)u_\la^\pm$
with unitaries $U_\la^\pm\in\Hom(\a^\pm_\la a,a\la)$
and $u_\la^\pm\in\Hom(\la\co a,\co a \a^\pm_\la)$.
Therefore we find
\[
\varphi_a \circ \a^\pm_{a;\la} \circ \varphi_a^{-1}
= {\co a}^{\,-1} \circ
\Ad(\co a(U_\la^\pm)u_\la^\pm) \circ \la \circ {\co a}
= \Ad (U_\la^\pm) \circ \a^\pm_\la ,
\]
and consequently maps
$\Hom(\a^+_\la,\a^\pm_\mu)\rightarrow
\Hom(\a^+_{a;\la},\a^\pm_{a;\mu})$,
$t\mapsto\varphi_a^{-1}(U_\mu^\pm t (U_\la^+)^*)$
are isomorphisms. In particular, the coupling matrix
arising from $N\subset M_a$ is the same as we obtained from
$N\subset M$. So here we have found some redundancy
for modular invariants from subfactors: Different
subfactors can produce the same coupling matrix $Z$,
and if we start with a given braided subfactor $N\subset M$,
than we obtain an irreducible subfactor for each
morphism in $\NXM$ producing the same $Z$,
though their Jones indices may well be different.
The simplest example is the trivial subfactor $N=M$,
where we obtain Jones extensions
$\la(N)\subset N\subset M_\la$ for each $\la\in\NXN$,
and they all will give us the trivial modular invariant
$Z_{\la,\mu}=\delta_{\la,\mu}$
(cf.\ \cite[Subsect.\ 6.2]{BEK2}).

\begin{figure}[htb]
\begin{center}
\unitlength 0.5mm
\begin{picture}(140,20)
\thinlines 
\multiput(0,10)(20,0){4}{\circle*{3}}
\multiput(100,10)(20,0){2}{\circle*{3}}
\thicklines
\multiput(0,10)(20,0){3}{\line(1,0){20}}
\put(60,10){\line(1,0){5}}
\put(95,10){\line(1,0){25}}
\dottedline{2}(65,10)(95,10)
\put(130,10){\circle{20}}
\put(0,4){\makebox(0,0){$1$}}
\put(20,4){\makebox(0,0){$2$}}
\put(40,4){\makebox(0,0){$3$}}
\put(60,4){\makebox(0,0){$4$}}
\put(100,4){\makebox(0,0){$\ell\,$--$1$}}
\put(117,4){\makebox(0,0){$\ell$}}
\end{picture}
\end{center}
\caption{Tadpole graph T$_\ell$, $\ell=2,3,4,...$}
\label{tadp}
\end{figure}
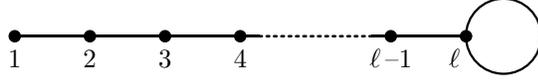

It is instructive to use these observations to demonstrate
that braided subfactors corresponding to the $\SUz_k$
systems can only produce modular invariant coupling matrices
with diagonal entries given as Coxeter exponents of
A-D-E Dynkin diagrams --- even if we would not know anything
about the list of modular invariants \cite{CIZ1,CIZ2,Kt}.
So suppose we have a given braided subfactor $N\subset M$
for the system $\NXN$ corresponding to the $\LSUz$ loop
group model at a level $k=1,2,3,...$, so that in
particular the endomorphisms are labelled by spins
$j=0,1,2,...,k$, we have fusion rules
\[
N_{j,j'}^{j''} = \left\{ \begin{array}{lc}
1 \qquad & |j-j'| \le j'' \le \min(j+j',2k-j-j') ,
\quad j+j'+j''\in2\bbZ , \\ 0 & {{\rm{otherwise}}},
\end{array} \right.
\]
and the statistics phases are given by
\[ \omega_j = \exp(2\pi\I h_j) , \qquad
h_j = \frac{j(j+2)}{4k+8} . \]
Now consider the (adjacency matrix of the)
$M$-$N$ fusion graph $G_1$,
i.e.\ the matrix of \erf{Gmatrix} corresponding to
the spin $j=1$, and let us wonder how it might look like.
Since we deal with the generator we will obtain a
connected graph, and since $\a$-induction preserves
statistical dimensions, we know already that its largest
eigenvalue only be $d_1=2\cos(\pi/(k+2))<2$.
This already forces
$G_1$ to be one of the A-D-E Dynkin diagrams with
dual Coxeter number $h=k+2$ or if $k=2\ell-1$ is
odd it could also be a tadpole T$_\ell$
which is displayed in \fig{tadp}
(see e.g.\ \cite{DiF} for such arguments).
The eigenvalues of T$_\ell$ are all simple and
given by $2\cos((r+1)\pi/(k+2))$ with
(``exponents'') $r=0,2,4,...,k-1$.
So if we pretend not to know that there is no modular
invariant at level $k=2\ell-1$ with $Z_{j,j}=1$
if $j$ is even and $Z_{j,j}=0$ if $j$ is odd,
then we have to rule out the possibility that there is
a subfactor with T$_\ell$ as $M$-$N$ fusion graph $G_1$.
This is easy: Assume that there is a subfactor $N\subset M$
which produces a coupling matrix with these diagonal entries.
Note that T$_\ell$ is A$_{2\ell}/\bbZ_2$ without a
fixed point, and thus the Perron-Frobenius weights
of each vertex of T$_\ell$ which is labelled by $j+1$
in \fig{tadp} is the same as for the unfolded A$_{2\ell}$,
i.e.\ $\sin((j+1)\pi/(k+2))/\sin(\pi/(k+2))$.
These numbers must be the statistical dimensions, up to an
overall normalizing factor which is fixed by the condition
that the global index of $\MXN$ is the same as $w$,
the global index of $\NXN$ (see e.g.\ \cite[p.\ 465]{BEK1}).
It is then easy to see that the statistical dimension of
the $M$-$N$ morphism which corresponds to the extremal
vertex labelled by ``$1$'' in \fig{tadp} is $\sqrt2$.
By the above arguments, there must hence be subfactor
which produces the same coupling matrix and which has
index two, i.e.\ is the unique $N\subset N\rtimes\bbZ_2$.
This would imply that there is an automorphism $\sigma\in\NXN$
such that $\sigma^2=\id$ and $[\sigma]\neq[\id]$.
But the only non-trivial automorphism in the $N$-$N$ system
is the spin $j=k$ simple current, however, by Rehren's lemma
\cite[Lemma 4.4.]{R3} this one cannot fulfill $\sigma^2=\id$
because its conformal dimension $k/4$ does not give a
statistics phase which is a second root of unity ---
contradiction.

\section{Structure of modular invariants from subfactors:\\
Type I coupling matrices and fusion rule isomorphisms}

Next we turn to the discussion of the distinction of
\typei\ and \typeii\ invariants in the subfactor framework.
In our general setting we have
\[
Z_{\la,\mu} = \langle\a^+_\la,\a^-_\mu\rangle
= \sum_{\tau\in\MXMo} b^+_{\tau,\la} b^-_{\tau,\mu} \,,
\]
with chiral branching coefficients
$b^\pm_{\tau,\la}=\langle\tau,\a^\pm_\la\rangle$.
To get this in the form of \erf{extaut} we would need
$b^-_{\tau,\la}=b^+_{\vartheta(\tau),\la}$ for
a permutation of the extended system, being identified
as the neutral system $\MXMo$.
Note that by $\vartheta(0)=0$ and $b^\pm_{\tau,0}=\delta_{\tau,0}$
(do not worry that we denote both the original and
the extended ``vacuum'' i.e.\ identity morphism
by the same symbol ``$0$'')
we are automatically forced to have symmetric vacuum coupling
$Z_{\la,0}=Z_{0,\la}$.
This corresponds basically to an ``extension of the chiral
algebras by primary fields'', namely those which appear in
the vacuum column or equivalently in the vacuum row of the
coupling matrix $Z$.
However, there are more general cases as we shall see,
which correspond to different extensions for the left
and right chiral algebra, and then the vacuum column of
the coupling matrix will be different from its vacuum row.
In this case we will be lead to different labelling sets of
extended sectors so that the extended modular invariant
coupling matrix is
\[
Z^\ext_{\tau_+,\tau_-} = \delta_{\tau_+,\vartheta(\tau_-)} \,,
\]
where $\vartheta$ is now an isomorphism between the two sets
of extended fusion rules, still subject to $\vartheta(0)=0$.
Note that when we have two different labelling sets it
makes no sense to ask whether a coupling matrix is
symmetric or not.

In the subfactor context, the crucial condition for \typei\ 
oupling matrices turns out to be the ``chiral locality condition'',
namely the following eigenvalue condition for
the statistics operator $\eps^+(\canr,\canr)$:
\begin{equation}
\label{chirloc}
\eps^+(\canr,\canr) \can(v) = \can(v)
\end{equation}
The name chiral locality was given since it was shown
in \cite{LR} that \erf{chirloc} is a necessary and sufficient
condition for locality of the extended net of (here: chiral)
observables in the nets of subfactors framework.
When chiral locality does hold then \cite[Prop.\ 3.3]{BE3}
\[
\langle\a^\pm_\la,\beta\rangle
=\langle\la,\sigma_\beta\rangle \,,
\]
whenever $\beta\in\MXMpm$. In particular, when
$\beta=\tau$ is neutral, i.e.\ lies in the intersection
$\MXMo=\MXMp\cap\MXMm$, then
\[
b^+_{\tau,\la}=b^-_{\tau,\la}\equiv b_{\tau,\la}\,,
\]
and we have a block decomposition or ``\typei''
modular invariant
\[
Z_{\la,\mu} = \sum_{\tau\in\MXMo} b_{\tau,\la}b_{\tau,\mu} \,.
\]
Now let us characterize to the other extreme case in the
subfactor context: pure permutation invariants.
As a corollary of \erf{chirglob}, the following conditions
are equivalent \cite[Prop.\ 3.2]{BEK2}:
\begin{itemize}
\item we have $Z_{\la,0}=\delta_{\la,0}$,
\item we have $Z_{0,\la}=\delta_{\la,0}$,
\item we have $\MXMo=\MXM$,
\item The coupling matrix $Z$ is a permutation, fixing the
  vacuum and corresponding to a fusion rule automorphism.
\end{itemize}
For the general case we would like to decompose a modular
invariant into its two parts, a \typei\ part together with a twist,
and in order to take care of heterotic vacuum coupling we
will need to implement such a twist by an isomorphism rather
than an automorphism. First we characterize chiral locality.
Indeed the following conditions are
equivalent \cite[Prop.\ 3.2]{BE4}:
\begin{itemize}
\item We have $Z_{\la,0}=\lan\canr,\la\ran$ for all $\la\in\NXN$.
\item We have $Z_{0,\la}=\lan\canr,\la\ran$ for all $\la\in\NXN$.
\item Chiral locality holds: $\eps^+(\canr,\canr)v^2=v^2$.
\end{itemize}
In other words: chiral locality holds if and only if the
dual canonical endomorphism is entirely
``visible'' in the vacuum row (and hence column)
of the coupling matrix.
Using results on intermediate subfactors from \cite{ILP},
we showed the following for a braided subfactor $N\subset M$
producing a coupling matrix $Z$ (no matter whether the
braiding is non-degenerate or not):
There are always \cite[Thm.\ 4.7]{BE4} intermediate subfactors
$N\subset M_\pm\subset M$, where $N\subset M_+$ and
$N\subset M_-$ fulfill the chiral locality condition
and produce coupling matrices (``\typei\ parents'')
$Z^+$ and $Z^-$, respectively, such that
$Z_{\la,0}=Z^+_{\la,0}=Z^+_{0,\la}$
$Z_{0,\la}=Z^-_{\la,0}=Z^-_{0,\la}$.
Moreover \cite[Thm.\ 5.5]{BE4}, the neutral systems
arising from $N\subset M$ and $N\subset M_\pm$ are canonically
isomorphic, and whenever $M_+=M_-$ then we have in fact
$Z^+_{\la,\mu}=\sum_\tau b_{\tau,\la} b_{\tau,\mu}=Z^-_{\la,\mu}$
and $Z_{\la,\mu}=\sum_\tau b_{\tau,\la} b_{\vartheta(\tau),\mu}$
is of Moore-Seiberg \cite{MS2}, Dijkgraaf-Verlinde \cite{DV} form,
where now the fusion rule automorphism arises from the
isomorphic neutral systems.
It is however important to notice that $M_+\neq M_-$,
even $Z^+\neq Z^-$ and $Z_{\la,0}\neq Z_{0,\la}$ can occur.
In fact we realized in \cite{BE4} the following coupling matrix
for $\mathit{SO}(16)_1$ (actually $\mathit{SO}(n)_1$ with
$n$ any multiple of $16$) from some subfactor $N\subset M$:
\[ 
Z = \left( \begin{array}{rrrr}
1 & 0 & 0 & 1 \\ 0 & 0 & 0 & 0 \\
1 & 0 & 0 & 1 \\ 0 & 0 & 0 & 0
\end{array}  \right) .
\]
The intermediate subfactors $N\subset M_+$ and
$N\subset M_-$ produce \typei\ parent coupling matrices
\[
Z^+ = \left( \begin{array}{rrrr}
1 & 0 & 1 & 0 \\ 0 & 0 & 0 & 0 \\
1 & 0 & 1 & 0 \\ 0 & 0 & 0 & 0
\end{array}  \right) \qquad \mbox{and} \qquad
Z^- = \left( \begin{array}{rrrr}
1 & 0 & 0 & 1 \\ 0 & 0 & 0 & 0 \\
0 & 0 & 0 & 0 \\ 1 & 0 & 0 & 1
\end{array}  \right), \]
respectively. (The first labels for rows and columns always
correspond to the vacuum.)
Because these invariants can be realized from subfactors
these coupling matrices are not spurious in the sense
that there are 2D RCFT's such that they encode the coupling
of left and right chiral sectors --- more about this in
the following section.

So what is the physical significance of the intermediate
factors $M_+$ and $M_-$? This can be understood if one
considers certain ``canonical tensor product subfactors''
$N\otimes N^\op\subset B$ which are directly related to
the possible existence of some 2D RCFT containing chiral
subtheories described by $N$ and encoded in the coupling
matrix $Z$ \cite{R5,R7}. Then $M_+\otimes M_-^\op\subset B$
turns out to be intermediate, and in the physical
interpretation of \cite{R5}, $M_+$ and $M_-$ correspond
precisely to the maximally extended chiral algebras
(in a sensible meaning). For more details, see \cite{BE4}.

\section{On the existence of 2D RCFT's and the
realization of modular invariants from subfactors}

The situation for modular invariants from subfactors can
simply be stated as follows:
For a given \typeiii\ von Neumann factor $N$ equipped with
a braided system of endomorphism $\NXN$, any embedding
$N\subset M$ of $N$ in a larger factor $M$ which is
compatible with the system $\NXN$ (in the sense that
the dual canonical endomorphism decomposes in $\NXN$)
defines a coupling matrix $Z$ through $\a$-induction.
This matrix $Z$ commutes with the matrices $Y$ and $\Omega$
arising from the braiding and in turn is a
``modular invariant'' whenever the braiding
is non-degenerate.
Suppose we start with a system corresponding to some
known RCFT data. More concretely, let us consider the
situation that our factor $N$ arises as a local factor
$N=N(I_\circ)$ of a  conformally (here: M\"obius) covariant
net $\{N(I)\}$ over $\bbR$ (or equivalently $S^1\setminus\zeta$),
as for example of the above sketched $\SUn_k$ loop group model.
Then the following questions are natural:
\begin{problem}
\label{p1}
Is any coupling matrix which can be produced by some brai\-ded
extension $N\subset M$ a physical invariant?
\end{problem}
And conversely:
\begin{problem}
\label{p2}
Can any physical modular invariant be
realized from some brai\-ded extension $N\subset M$?
\end{problem}
The first difficulty here is that one needs to specify what
the term ``physical'' means. Quite often in the literature,
any matrix commuting with $S$ and $T$ and subject
to the constraint that all entries are non-negative integers
and with normalization $Z_{0,0}=1$ is called a physical invariant.
Well, with this interpretation of ``physical'' the solution of
Problem \nolinebreak \ref{p1} is trivially the answer ``Yes''
since we have already established that our coupling matrices
have these properties.
It is also not too difficult to see that the solution of
Problem \nolinebreak \ref{p2} is just ``No'' with this
interpretation of ``physical'': Namely, our
general theory says that there is always some associate
extended theory carrying another representation of
the modular group $\SLZ$ which is compatible with
the chiral branching rules (see \cite[Sect.\ 6]{BE4}).
It is however known \cite{SY2,V,FSS}
that there are ``spurious'' modular invariants satisfying
the above constraints but which do not admit an extended
modular S-matrix.

Another, physically much more interesting specification of
``physical'' (but unfortunately mathematically harder to reach)
is that $Z$ arises from ``the existence of some 2D RCFT''.
A reasonable way of making this precise seems for us to
be the concept of chiral observables as light-cone nets
built in an observable net over 2D Minkowski space \cite{R5}.
And in fact Rehren has shown \cite{R7} that in the above
situation any braided extension $N\subset M$
determines an entire local 2D conformal field theory over
Minkowski space, and that indeed the vacuum Hilbert space
of the 2D net decomposes upon restriction to the tensor
product of the left and right chiral observables according
to \erf{statedecom} with $Z$ being precisely the
matrix arising from $N\subset M$ through $\a$-induction.
Therefore we obtain a positive solution for
Problem \nolinebreak \ref{p1}
even with this more subtle notion of a ``physical invariant''.

Now let us turn to the converse direction,
Problem \nolinebreak \ref{p2}.
It should be noticed that there can be local 2D extensions
of tensored left and right chiral observables
(besides the trivial extension) which are completely
compatible with conformal symmetry and whose coupling
matrices do commute with $T$ (due to locality) but are
not S-invariant (see Rehren's contribution to this volume).
However, Problem \nolinebreak \ref{p2} is restricted to
modular invariant coupling matrices, and considering
T- and S-invariant 2D extensions,
we tend to believe that the solution to
Problem \nolinebreak \ref{p2} is again a ``Yes''.
However, this requires a proof!
But since the general classification problem for
modular invariants is not solved in general and since
it is still a quite subtle question to distinguish
physical from spurious invariants even if somebody
provides you with a complete list of normalized,
non-negative integer matrices in the commutant of $S$ and $T$
for some given model, for the time being the state of the art
seems to admit only the following recipe:
\begin{enumerate}
\item Pick the first $Z$ from the list and
 \begin{enumerate}
 \item either realize it from a subfactor (and possibly
       classify inequivalent realizations),
 \item or disprove the existence of a subfactor producing
       this coupling matrix.
 \end{enumerate}
\item Pick the next $Z$.
\item If you are lucky, do these steps for entire classes
 of $Z$'s rather than for single ones, so that you can cover
 certain parts of any list.
\end{enumerate}
Tackling step 1.(a) with a case-by-case analysis has been
carried out for a few models. A relation between the
A-D-E modular invariants of $\SUz$ and A-D-E
Goodman-de la Harpe-Jones subfactors \cite{GHJ}
can be found in \cite{O}. The \typeiii\ subfactors
for E$_6$, E$_8$ and $D_4$ where analyzed in \cite{X1}
and found to produce these Dynkin diagrams as fusion
diagrams of sectors. The subfactors for the entire
D$_{\mathrm{even}}$ series were constructed in \cite{BE2}
and shown to produce the diagrams as fusion graphs,
and that the modular invariants can be recovered by the
formula \erf{Zdef} was proven in \cite{BE3}.
The remaining cases were treated in the $\a$-induction
setting extensively in \cite{BEK2}.
For $\SUd$, the \typei\ exceptional cases were first
analysed in the subfactor context in \cite{X1},
further analysis also covering the $\cD$ series was
carried out in \cite{BE2,BE3,BEK2}, and
Ocneanu claimed the solution of the $\SUd$
problem in January 2000.

Anyway, let us try to tackle step 1.(a) in the more efficient
way, namely looking at classes of modular invariants.
Our discussion will be focussed on $\SUn_k$,
but the general arguments can also be translated to other models.
First of all the trivial invariants,
$Z_{\la,\mu}=\delta_{\la,\mu}$, are obtained from
the trivial subfactor $N\subset M$ with $M=N$.
Next, there are the exceptional modular invariants
arising from conformal inclusions.
A conformal inclusion means that the level $1$ representations
of some loop group of a Lie group restrict in a finite manner
to the positive energy representations of a certain embedded
loop group of an embedded (simple) Lie group at some level.
For $\SUz$, the modular invariants arising from conformal
embeddings are the E$_6$, E$_8$ and D$_4$ ones,
corresponding to embeddings $\SUz_{10}\subset\SOf_1$,
$\SUz_{28}\subset(\Gtwo)_1$ and $\SUz_4\subset\SUd_1$, respectively.
(The latter happens to be a simple current invariant at the same time.)
For $\SUd$, the invariants from conformal embeddings
are $\cD^{(6)}$, $\cE^{(8)}$, $\cE^{(12)}$ and $\cE^{(24)}$,
corresponding to $\SUd_3\subset{\mathit{SO}}(8)_1$,
$\SUd_5\subset{\mathit{SU}}(6)_1$, $\SUd_9\subset(\rmE_6)_1$,
$\SUd_{21}\subset(\rmE_7)_1$, respectively.
By taking such an embedding as a local subfactor in the vacuum
representation, any conformal inclusion determines a braided
subfactor of finite index (see \cite{Wup,R4,RST,LR,X1,BE2}), 
which in turn produces a modular invariant, being precisely
the \typei\ (since the embedding level one theory is always
local) exceptional invariant which arises from the
diagonal invariant of the extended theory \cite{BE3}.
So here a class of exceptional modular invariants
is covered in the subfactor context at one stroke,
and the consequently existing extended RCFT is
of course the level \nolinebreak $1$ theory of the larger
affine Lie algebra.

The situation is even better for simple current invariants,
which in a sense produce the majority of non-trivial
modular invariants. Simple currents \cite{SY1} are primary fields
with unit quantum dimension and appear in our framework
automorphisms in the system $\NXN$. They form a closed abelian
group $G$ under fusion which is hence a product of
cyclic groups. Simple currents give rise to modular invariants,
and all such invariants have been classified \cite{GRS,KS}.
As focussing on $\SUn$ here, we will simply consider
cyclic simple current groups $\bbZ_n$.

By taking a generator $[\sigma]$ for of the cyclic
simple current group $\bbZ_n$ we can construct the
crossed product subfactor $N\subset M=N\rtimes\bbZ_n$
whenever we can choose a representative $\sigma$
in each such simple current sector such that we have
exact cyclicity $\sigma^n=\id$ (and not only as sectors).
As we are starting with a chiral quantum field theory,
Rehren's lemma \cite{R3} applies which states that such a
choice is possible if and only if the statistics phase
$\om_\sigma$ is an $n$-th root of unity, i.e.\ if and
only if the conformal weight $h_\sigma$ is an integer
multiple of $1/n$.
Sometimes this may only be possible for a simple current
subgroup $\bbZ_m\subset\bbZ_n$ (with $m$ a divisor of $n$)
but any such non-trivial ($m\neq 1$) subgroup
gives rise to a non-trivial subfactor and in turn to
a modular invariant.
It is easy to see that in fact all simple current invariants
are realized this way. For $\SUn_k$ the simple current
group $\bbZ_n$ corresponds to weights $k\Lambda_{(j)}$,
$j=0,1,...,n-1$.
The conformal dimensions are $h_{k\Lambda_{(j)}}=kj(n-j)/2n$,
which by Rehren's Lemma \cite{R3} allow for full $\bbZ_n$
extensions except when $n$ is even and $k$ is odd in which
case the maximal extension is $N\subset M=N\rtimes\bbZ_{n/2}$
because we can only use the even labels $j$.
(This reflects the fact that e.g.\
for $\SUz$ there are no D-invariants at odd levels.)
Thus Rehren's lemma has told us that extensions are
labelled by all the divisors of $n$ unless $n$ is even and $k$
is odd in which case they are labelled by the divisors of $n/2$.
This matches exactly the simple current modular invariant
classification of \cite{GRS,KS}.
An extension by a simple current subgroup $\bbZ_m$,
with $m$ is a divisor of $n$ or $n/2$, is moreover local, if the
generating current (and hence all in the $\bbZ_m$ subgroup)
has integer conformal weight, $h_{k\Lambda_{(q)}}\in\bbZ$,
where $n=mq$. This happens exactly if $kq\in 2m\bbZ$ if $n$
is even, or $kq\in m\bbZ$ if $n$ is odd \cite{BE3}.
For $\SUz$ this corresponds to the $\rmD_{\mathrm{even}}$
series whereas the $\rmD_{\mathrm{odd}}$
series are non-local extensions. For $\SUd$, there is a simple
current extension at each level, but only those at $k\in 3\bbZ$
are local. Clearly, the cases with chiral locality match
exactly the \typei\ simple current modular invariants.

With these techniques we can obtain a large number of
modular invariants from subfactors. Nevertheless we still
do not have a systematic procedure to get \textit{all} physical
invariants. The more problematic cases are typically the
exceptional \typeii\ invariants. Therefore let us now
tackle a large class of exceptional \typeii\ invariants,
namely those which are \typeii\ descendants of conformal
embeddings.

\section{Type II descendants of modular
 invariants from conformal inclusions}

Since for any conformal inclusion the ambient theory is
described by the level $1$ representations of the
embedding loop group and therefore is typically a pure
simple current theory (whenever simply laced Lie groups
are dealt with), the subfactors producing their modular
invariants can be constructed by simple current methods.
Therefore we will obtain the relevant subfactors for
\typeii\ decendant modular invariants, e.g.\ the
conjugation $CZ$ of a conformal inclusion invariant $Z$
through crossed products.

For a while we will be looking at the so-called $\bbZ_n$
conformal field theories as treated in \cite{D}, which have
$n$ sectors, labelled by $j=0,1,2,...,n-1$ (mod$\,n$),
obeying $\bbZ_n$ fusion rules, and conformal dimensions
of the form $h_j=a j^2/2n$ (mod$\,1$), where $a$ is an integer
mod$\,2n$, $a$ and $n$ coprime and $a$ is even whenever
$n$ is odd. The modular invariants of such models have been
classified \cite{D}. They are labelled by the divisors
$\delta$ of $\tn$, where $\tn=n$ if $n$ is odd and $\tn=n/2$ if
$n$ is even. Explicitly, the modular invariants $Z^{(\delta)}$
are given by
\[
Z^{(\delta)}_{j,j'} = \left\{ \begin{array}{ll}
1 & \qquad \mbox{if}\;\;\; j,j'=0\,\mbox{mod}\,\alpha
\;\;\;\mbox{and}\;\;\; j'=\omega(\delta)j \,\mbox{mod}\,
n/\alpha\,, \\[.4em]
0 & \qquad \mbox{otherwise}\,, \end{array} \right.
\]
where $\a=\mathrm{gcd}(\delta,\tn/\delta)$ so that there
are numbers $r,s\in\bbZ$ such that
$r\tn/\delta\a-s\delta/\a=1$ and then $\omega(\delta)$
is defined as $\omega(\delta)=r\tn/\delta\a+s\delta/\a$.
The trivial invariant corresponds to $\delta=\tn$,
i.e.\ $Z^{(\tn)}=\bfe$ and $\delta=1$ gives the charge
conjugation matrix, $Z^{(1)}=C$.

It is straightforward combinatorics \cite{BE5} to show that
\[
Z^{(\delta)}_{j,j} = \left\{ \begin{array}{ll}
1 & \qquad \mbox{if}\;\;\; j=0\,\mbox{mod}\,\tn/\delta, \\[.4em]
0 & \qquad \mbox{otherwise}\,. \end{array} \right.
\]
This yields $\tr(Z^{(\delta)}) = \epsilon\delta$
for the trace of $Z^{(\delta)}$; here $\epsilon=2$ if
$n$ is even and $\epsilon=1$ if $n$ is odd.
Now suppose that for such a $\bbZ_n$ theory at hand
we have corresponding braided automorphisms $\tau_j$
of some \typeiii\ factor $N$, obeying $\bbZ_n$ fusion rules
and such that their statistical phases are given by
$\E^{2\pi\I h_j}$ with conformal weights $h_j$ as above
(as is the case e.g.\ for $\SUn$ level 1 theories).
Note that if $n$ is odd then we can always assume that
$\tau_1^n=\id$ as morphisms (and our system can be
chosen as $\{\tau_1^j\}_{j=0}^{n-1}$). However,
if $n$ is even, then we cannot choose a representative
of the sector $[\tau_1]$ such that its $n$-th power gives
the identity, nevertheless we can always assume that
$\tau_\epsilon^\tn=\id$.
Thus we have a simple current (sub-) group $\bbZ_\tn$,
for which we can form the crossed product subfactor
$N\subset M=N\rtimes\bbZ_{\tn/\delta}$ for any divisor
$\delta$ of $\tn$. It is quite easy to see that
$N\subset M=N\rtimes_{\tau_{\epsilon\delta}}\bbZ_{\tn/\delta}$
indeed realizes $Z^{(\delta)}$: The crossed product by
$\bbZ_{\tn/\delta}$ gives the dual canonical endomorphism sector
$[\canr]=[\id]\oplus[\tau_{\epsilon\delta}]\oplus
[\tau_{\epsilon\delta}^2]\oplus\ldots\oplus
[\tau_{\epsilon\delta}^{\tn/\delta-1}]$.
The formula
$\langle\iota\tau_j,\iota\tau_{j'}\rangle=
\langle\canr\tau_j,\tau_{j'}\rangle$ then shows that
the system of $M$-$N$ morphisms is labelled by
$\bbZ_n/\bbZ_{\tn/\delta}\simeq\bbZ_{\epsilon\delta}$,
i.e.\ $\#\MXN=\epsilon\delta$.
Therefore our general theory implies that the modular invariant
arising from $N\subset M=N\rtimes\bbZ_{\tn/\delta}$ has trace
equal to $\epsilon\delta$, and thus must be $Z^{(\delta)}$.
Consequently all modular invariants classified in \cite{D}
are realized from subfactors.

It is instructive to apply the above results to descendant
modular invariants of conformal inclusions. Let us consider
the conformal inclusion $\SUf_6\subset{\mathit{SU}}(10)_1$.
The associated modular invariant, which can be found
in \cite{SY1}, reads
\[
\cZ = \sum_{j\in\bbZ_{10}} |\chi^j|^2
\]
with ${\mathit{SU}}(10)_1$ characters decomposing into
$\SUf_6$ characters as
\[
\begin{split}
\chi^0 = \chi_{0,0,0} + \chi_{0,6,0} + \chi_{2,0,2}
 + \chi_{2,2,2}, \;\; &
\chi^5 = \chi_{0,0,6} + \chi_{6,0,0} + \chi_{0,2,2} + \chi_{2,2,0},\\
\chi^1 = \chi_{0,0,2} + \chi_{2,4,0} + \chi_{2,1,2} , \quad &
\chi^6 = \chi_{4,0,0} + \chi_{0,2,4} + \chi_{1,2,1} ,\\
\chi^2 = \chi_{0,1,2} + \chi_{2,3,0} + \chi_{3,0,3} , \quad &
\chi^7 = \chi_{3,0,1} + \chi_{1,2,3} + \chi_{0,3,0} ,\\
\chi^3 = \chi_{1,0,3} + \chi_{3,2,1} + \chi_{0,3,0} , \quad &
\chi^8 = \chi_{0,3,2} + \chi_{2,1,0} + \chi_{3,0,3} ,\\
\chi^4 = \chi_{0,0,4} + \chi_{4,2,0} + \chi_{1,2,1} , \quad &
\chi^9 = \chi_{2,0,0} + \chi_{0,4,2} + \chi_{2,1,2} .
\end{split}
\]
We observe that $Z$ has 32 diagonal entries.
As usual, this invariant can be realized from the conformal
inclusion subfactor
\[
\pi^0({\mathit{L}}_I{\mathit{SU}}(4))'' \subset
\pi^0({\mathit{L}}_I{\mathit{SU}}(10))'' ,
\]
with $\pi^0$ denoting the level 1 vacuum representation
of $\mathit{LSU}(10)$. We will denote this subfactor by
$N\subset M_+$. The dual canonical endomorphism
sector corresponds to the vacuum block:
\[
[\canr_+] = [\la_{0,0,0}]\oplus[\la_{0,6,0}]\oplus
[\la_{2,0,2}]\oplus[\la_{2,2,2}] .
\]
Proceeding with $\a$-induction
$\la_{p,q,r}\mapsto\a^\pm_{+;p,q,r}\in\End(M_+)$,
it is a straightforward calculation that the graphs
describing left multiplication by fundamental generators
$[\a^\pm_{+;1,0,0}]$ and $[\a^\pm_{+;0,1,0}]$ (which is
the same as right multiplication by $[\la_{1,0,0}]$ and
$[\la_{0,1,0}]$, respectively) on the system of $M_+$-$N$
sectors gives precisely the graphs found by Petkova and
Zuber \cite[Figs.\ 1 and 2]{PZp} by their more empirical
procedure to obtain graphs with spectrum matching the diagonal
part of some given modular invariant. In our framework,
the graph \cite[Fig.\ 1]{PZp} obtains the following meaning:
Take the outer wreath, pick a vertex with 4-ality 0 and
label it by $[\iota_+]\equiv[\tau_0\iota_+]$,
where $\iota_+:N\hookrightarrow M_+$ denotes the
injection homomorphism, as usual.
Going around in a counter-clockwise direction the vertices
will then be the marked vertices labelled by the $\bbZ_{10}$
sectors $[\tau_1\iota_+]$, $[\tau_2\iota_+]$, .... ,
$[\tau_9\iota_+]$ of $\mathit{SU}(10)_1$.
Passing to the next inner wreath
the 4-ality 1 vertex adjacent to $[\iota_+]$ is then
the sector $[\a^\pm_{+;1,0,0}\iota_+]=[\iota_+\la_{1,0,0}]$,
and the others its $\bbZ_{10}$ translates.
Similarly the inner wreath consists
of the $\bbZ_{10}$ translates of $[\iota_+\la_{0,1,0}]$.
The remaining two vertices in the center correspond to
subsectors of the reducible $[\iota\la_{1,1,0}]$ and
$[\iota\la_{0,1,1}]$. The graph itself then represents
left (right) multiplication by $[\a^\pm_{+;1,0,0}]$
($[\la_{1,0,0}]$).

As for $\mathit{LSU}(10)$ at level 1 we are in fact dealing
with a $\bbZ_n$ conformal field theory, we have $n=10$ and
$\tn=5$, thus we know that there are only two modular
invariants: The diagonal one which in restriction to
$\mathit{LSU}(4)$ gives exactly the above \typei\ invariant
$Z\equiv Z^{(5)}$, but there is also the charge conjugation
invariant $CZ\equiv Z^{(1)}$, written in characters as
\[
C\cZ = \sum_{j\in\bbZ_{10}} \chi^j (\chi^{-j})^* .
\]
Whereas $Z^{(5)}$ can be thought of as the trivial
extension $M_+\subset M_+$, the conjugation invariant $Z^{(1)}$
can be realized from the crossed product
$M_+\subset M=M_+\rtimes\bbZ_5$ which has dual canonical
endomorphism sector
\[
[\canr^\ext]=[\tau_0]\oplus[\tau_2]\oplus[\tau_4]
\oplus[\tau_6]\oplus[\tau_8] .
\]
So far we have considered the situation on the ``extended level'',
but we may now descend to the level of $\SUf_6$ sectors and
characters. Namely we may consider the subfactor
$N\subset M=M_+\rtimes\bbZ_5$. Its dual canonical endomorphism
sector $[\canr]$ is obtained by $\sigma$-restriction of $[\canr^\ext]$
which can now be read off from the character decomposition,
\[
\begin{split}
[\canr]&=\bigoplus_{j=0}^4 [\sigma_{\tau_{2j}}]
=[\la_{0,0,0}]\oplus[\la_{0,6,0}]\oplus[\la_{2,0,2}]
\oplus[\la_{2,2,2}]\oplus [\la_{0,1,2}]
\oplus[\la_{2,3,0}]\\[.4em]
& \quad
\oplus 2[\la_{3,0,3}]\oplus[\la_{0,0,4}]
\oplus[\la_{4,2,0}]\oplus 2[\la_{1,2,1}]
\oplus[\la_{4,0,0}]\oplus[\la_{0,2,4}]
\oplus[\la_{0,3,2}]\oplus[\la_{2,1,0}].
\end{split}
\]
This subfactor produces the conjugation invariant
$CZ$ written in $\SUf_6$ characters which is
the same as taking the original $\SUf_6$ conformal
inclusion invariant and conjugating on the level of
the $\SUf_6$ characters. Note that this invariant has
only 16 diagonal entries.

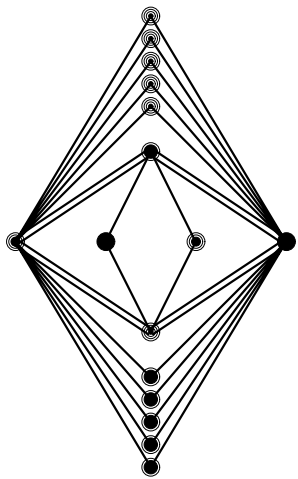
\begin{figure}[htb]
\begin{center}
\unitlength 0.6mm
\begin{picture}(80,110)
\thinlines
\multiput(30,50)(40,0){2}{\circle*{4}}
\multiput(10,50)(40,0){2}{\circle{4}}
\multiput(10,50)(40,0){2}{\circle{3}}
\multiput(10,50)(40,0){2}{\circle*{2}}
\put(40,70){\circle{4}}
\put(40,70){\circle*{3}}
\put(40,30){\circle{4}}
\put(40,30){\circle{3}}
\put(40,30){\circle{2}}
\put(40,30){\circle*{1}}
\multiput(40,80)(0,5){5}{\circle{4}}
\multiput(40,80)(0,5){5}{\circle{3}}
\multiput(40,80)(0,5){5}{\circle{2}}
\multiput(40,80)(0,5){5}{\circle*{1}}
\multiput(40,0)(0,5){5}{\circle{4}}
\multiput(40,0)(0,5){5}{\circle*{3}}
\thicklines
\path(30,50)(40,70)(50,50)(40,30)(30,50)
\path(10,50)(40,80)(70,50)(40,20)(10,50)
(40,85)(70,50)(40,15)(10,50)
(40,90)(70,50)(40,10)(10,50)
(40,95)(70,50)(40,5)(10,50)
(40,100)(70,50)(40,0)(10,50)
\path(9.5,50.5)(39.5,70.5)
\path(10.5,49.5)(40.5,69.5)
\path(9.5,49.5)(39.5,29.5)
\path(10.5,50.5)(40.5,30.5)
\path(39.5,30.5)(69.5,50.5)
\path(40.5,29.5)(70.5,49.5)
\path(39.5,69.5)(69.5,49.5)
\path(40.5,70.5)(70.5,50.5)
\end{picture}
\end{center}
\caption{Graph $G_1$ associated to the conjugation invariant
of the conformal inclusion $\SUf_6\subset{\mathit{SU}}(10)_1$}
\label{orbPZ1}
\end{figure}

When passing from $M_+$ to $M=M_+\rtimes\bbZ_5$,
the $M_+$-$N$ system will change
to the $M$-$N$ system in such a way
that all sectors which are translates
by $\tau_{2j}$, $j=0,1,2,3,4$, have to be identified,
and similarly fixed points split. Thus our new system of
$M$-$N$ morphisms will be some kind of orbifold of the old one.
To see this, we first recall that all the irreducible
$M_+$-$N$ morphisms are of the form $\beta\iota_+$ with
$\beta\in\MXMppm$. To such an irreducible $M_+$-$N$ morphism
$\beta\iota_+$ we can now associate an $M$-$N$ morphism
$\iota^\ext\beta\iota_+$ which may no longer be irreducible;
here $\iota^\ext$ is the injection homomorphism
$M_+\hookrightarrow M$. Then the reducibility can be
controlled by Frobenius reciprocity as we have
\[
\langle\iota^\ext\beta\iota_+,\iota^\ext\beta'\iota_+\rangle=
\langle\canr^\ext\beta\iota_+,\beta'\iota_+\rangle ,
\]
and $\canr^\ext=\co\iota ^\ext\iota^\ext$.
Carrying out the entire computation we find that
in contrast to the original 32 $M_+$-$N$ sectors
we are left with only 16 $M$-$N$ sectors, and the
right multiplication by $[\la_{1,0,0}]$ is displayed
graphically as in \fig{orbPZ1}.
Here the 4-alities 0,1,2,3 of the vertices are indicated
by solid circles of decreasing size. The $[\iota]$
vertex (with $\iota=\iota^\ext\iota_+$ denoting the
injection homomorphism $N\hookrightarrow M$ of the
total subfactor $N\subset M=M_+\rtimes \bbZ_5$)
is the 4-ality 0 vertex in the center of the picture,
and the 4-ality 1 vertex above corresponds to
$[\iota\la_{1,0,0}]$. Each group of five vertices on
the top and the bottom of the picture arise from the
splitting of the two central vertices of the graphs
in \cite{PZp} as they are $\bbZ_5$ fixed points.
That our orbifold graph inherits the 4-ality of the
original graph is due to the fact that all entries in
$[\canr]$ have 4-ality zero which in turn comes from the
fact that all even marked vertices (corresponding to
the subgroup $\bbZ_5\subset\bbZ_{10}$) of the graph of
Petkova and Zuber have 4-ality zero.
We also display the graph corresponding to the second
fundamental representation, namely the right multiplication
by $[\la_{0,1,0}]$ in \fig{orbPZ1}.
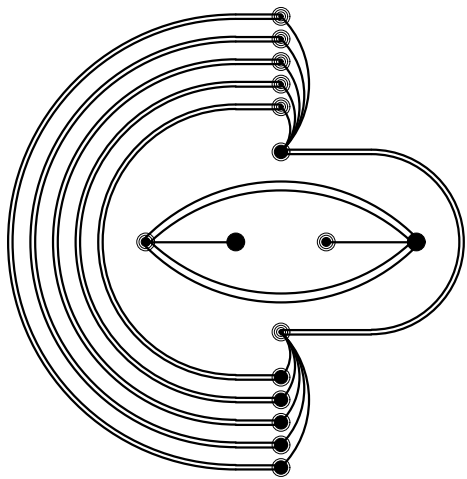
\begin{figure}[htb]
\begin{center}
\unitlength 0.6mm
\begin{picture}(120,110)
\thinlines
\multiput(60,50)(40,0){2}{\circle*{4}}
\multiput(40,50)(40,0){2}{\circle{4}}
\multiput(40,50)(40,0){2}{\circle{3}}
\multiput(40,50)(40,0){2}{\circle*{2}}
\put(70,70){\circle{4}}
\put(70,70){\circle*{3}}
\put(70,30){\circle{4}}
\put(70,30){\circle{3}}
\put(70,30){\circle{2}}
\put(70,30){\circle*{1}}
\multiput(70,80)(0,5){5}{\circle{4}}
\multiput(70,80)(0,5){5}{\circle{3}}
\multiput(70,80)(0,5){5}{\circle{2}}
\multiput(70,80)(0,5){5}{\circle*{1}}
\multiput(70,0)(0,5){5}{\circle{4}}
\multiput(70,0)(0,5){5}{\circle*{3}}
\thicklines
\path(40,50)(60,50)
\path(80,50)(100,50)
\put(70,21){\arc{84.853}{3.927}{5.498}}
\put(70,19){\arc{84.853}{3.927}{5.498}}
\put(70,81){\arc{84.853}{0.785}{2.356}}
\put(70,79){\arc{84.853}{0.785}{2.356}}
\put(65,25){\arc{14.142}{5.498}{0.785}}
\put(62.5,22.5){\arc{21.213}{5.498}{0.785}}
\put(60,20){\arc{28.284}{5.498}{0.785}}
\put(57.5,17.5){\arc{35.355}{5.498}{0.785}}
\put(55,15){\arc{42.426}{5.498}{0.785}}
\put(65,75){\arc{14.142}{5.498}{0.785}}
\put(62.5,77.5){\arc{21.213}{5.498}{0.785}}
\put(60,80){\arc{28.284}{5.498}{0.785}}
\put(57.5,82.5){\arc{35.355}{5.498}{0.785}}
\put(55,85){\arc{42.426}{5.498}{0.785}}
\put(60,50){\arc{61}{1.571}{4.712}}
\put(60,50){\arc{59}{1.571}{4.712}}
\put(60,50){\arc{71}{1.571}{4.712}}
\put(60,50){\arc{69}{1.571}{4.712}}
\put(60,50){\arc{81}{1.571}{4.712}}
\put(60,50){\arc{79}{1.571}{4.712}}
\put(60,50){\arc{91}{1.571}{4.712}}
\put(60,50){\arc{89}{1.571}{4.712}}
\put(60,50){\arc{101}{1.571}{4.712}}
\put(60,50){\arc{99}{1.571}{4.712}}
\path(60,0.5)(70,0.5)
\path(60,-0.5)(70,-0.5)
\path(60,5.5)(70,5.5)
\path(60,4.5)(70,4.5)
\path(60,10.5)(70,10.5)
\path(60,9.5)(70,9.5)
\path(60,15.5)(70,15.5)
\path(60,14.5)(70,14.5)
\path(60,20.5)(70,20.5)
\path(60,19.5)(70,19.5)
\path(60,80.5)(70,80.5)
\path(60,79.5)(70,79.5)
\path(60,85.5)(70,85.5)
\path(60,84.5)(70,84.5)
\path(60,90.5)(70,90.5)
\path(60,89.5)(70,89.5)
\path(60,95.5)(70,95.5)
\path(60,94.5)(70,94.5)
\path(60,100.5)(70,100.5)
\path(60,99.5)(70,99.5)
\put(90,50){\arc{41}{4.712}{1.571}}
\put(90,50){\arc{39}{4.712}{1.571}}
\path(70,30.5)(90,30.5)
\path(70,29.5)(90,29.5)
\path(70,70.5)(90,70.5)
\path(70,69.5)(90,69.5)
\end{picture}
\end{center}
\caption{Graph $G_2$ associated to the conjugation invariant
of the conformal inclusion $\SUf_6\subset{\mathit{SU}}(10)_1$}
\label{orbPZ2}
\end{figure}

Let us finally remark that the conformal inclusion
invariant $Z$ has the funny property $Z^*Z=3Z+CZ$.
This is remarkable as this is the first \typei\
invariant we have encountered which does not fulfill
$Z^*Z=x_+Z$, where $x_+=\sum_\la Z_{\la,0}^2$.
Since the diagonal part of $Z^*Z$ describes the spectral
properties of the fusion graphs of chiral generators
in the full system \cite{BEK1,BEK2} we expect that
for the conformal inclusion subfactor $N\subset M_+$
the fusion graph of the chiral generator $[\a^+_{1,0,0}]$
(or $[\a^-_{1,0,0}]$) in the full $M_+$-$M_+$ system has four
connected components (each of which corresponds to a nimrep),
three of them being the graph of \cite[Fig.\ 1]{PZp} and
one component being the graph in \fig{orbPZ1}.
In fact, it is easy to see that for any (non-degenerately)
braided subfactor $N\subset M$ the number of $\MXMo$
fusion orbits in $\MXMpm$ is equal to the number of
$\MXMmp$ fusion orbits in $\MXM$. Moreover, a simple
Perron-Frobenius argument shows that this number is
exactly $x_\pm$ (with $x_-=\sum_\la Z_{0,\la}^2$).

The conformal inclusion $\SUd_5\subset{\mathit{SU}}(6)_1$
can be treated along the same lines \cite{BE5}.
The associated $\SUd_5$ modular invariant, 
is labelled by the graph $\cE^{(8)}$.
The ambient ${\mathit{SU}}(6)_1$ has besides the diagonal
only the conjugation invariant which is the obtained
through a $\bbZ_3$ extension on top of the conformal
inclusion subfactor, and the $\bbZ_3$ quotient
collapses the 12 vertices of $\cE^{(8)}$ to 4,
yielding exactly the graph ${\cE^{(8)}}{}^*$ in the list
of Di Francesco and Zuber (see e.g.\ \cite{BPPZ}).
So with this procedure we understand quite
generally why the descendants of modular
invariants of conformal inclusions (whenever the ambient
theory has $\bbZ_n$ fusion rules) are in fact labelled
by orbifold graphs of the graph labelling the original,
block-diagonal conformal inclusion invariant, and
why the conjugation invariant corresponds to the
maximal $\bbZ_\tn$ orbifold.

In the above examples, the trivial and conjugation invariant
of the extended theory still remain distinct when written
in terms of the $\SUf_6$ characters. This need not be the
case in general. Let us look at a familiar modular invariant
of $\SUd$ at level 9, namely
\[
\cZ_{\cE^{(12)}} = \abs{ \chi_{0,0} + \chi_{9,0} +
\chi_{0,9} + \chi_{4,1} + \chi_{1,4} + \chi_{ 4,4} }^2
+ 2 \abs{ \chi_{2,2} + \chi_{5,2} + \chi_{2,5} }^2 ,
\]
which arises from the conformal embedding
$\SUd_9\subset(\rmE_6)_1$. Now $\rmE_6$
at level 1 gives a $\bbZ_3$ theory and in terms of
the extended characters the above invariant is the
trivial extended invariant
\[
\cZ_{\cE_1^{(12)}} = \abs{\chi^0}^2+\abs{\chi^1}^2+\abs{\chi^2}^2 ,
\]
using obvious notation. Here both the $(\rmE_6)_1$
characters $\chi^1$ and $\chi^2$ specialize to
$\chi_{2,2} + \chi_{5,2} + \chi_{2,5}$ in terms of
$\SUd_9$ variables. Let $N\subset M_+$ denote the conformal
inclusion subfactor obtained by analogous means as in the
previous example.\footnote{The notation $N\subset M_+$ for the
conformal inclusion subfactor which indicates that it will be
the maximal local intermediate subfactor (\`a la \cite{BE4})
of some extension $N\subset M=M_+\rtimes \bbZ_\ell$ using
simple currents is appropriate for the examples discussed
here but not in general. Other conformal inclusions as
for instance ${\mathit{SU}}(7)_7\subset{\mathit{SO}}(48)_1$
or ${\mathit{SU}}(8)_{10}\subset{\mathit{SU}}(36)_1$
can also have \typei\ descendants as coming from
\textit{local} simple current extensions of the ambient theory
(e.g.\ ${\mathit{SU}}(36)_1\rtimes\bbZ_3$ for the latter example).
In other words, for such descendant invariants,
the ambient affine Lie algebra does \textit{not}
provide the maximally extended chiral algebra which
then is actually larger.}
It has been treated in \cite{X1,BE3}
and produces the graph $\cE^{(12)}_1$ of the list of
Di Francesco and Zuber as chiral fusion graphs ---
and in turn as $M_+$-$N$ fusion graph,
thanks to chiral locality.

Corresponding to the two divisors
3 and 1 of 3, we know that besides the trivial there
is only the conjugation invariant of our $\bbZ_3$ theory.
It is given as
\[
C\cZ_{\cE_2^{(12)}} = \abs{\chi^0}^2+\chi^1(\chi^2)^*
+\chi^2(\chi^1)^*
\]
but this distinct invariant restricts to the same
invariant $Z_{\cE^{(12)}}$ when specialized to
$\SUd_9$ variables. Nevertheless we will obtain a different
subfactor $N\subset M$ since the conjugation invariant
of our $\bbZ_3$ theory is realized from the extension
$M_+\subset M=M_+\rtimes\bbZ_3$. In particular,
the subfactor $N\subset M$ has dual canonical
endomorphism sector
\[
[\canr] = [\la_{0,0}] \oplus [\la_{9,0}] \oplus
[\la_{0,9}] \oplus [\la_{4,1}] \oplus [\la_{1,4}] \oplus
[\la_{ 4,4}] \oplus 2 [\la_{2,2}] \oplus 2 [\la_{5,2}]
\oplus 2[\la_{2,5}] ,
\]
determined by $\sigma$-restriction of
\[
[\canr^\ext] = [\tau_0] \oplus [\tau_1] \oplus [\tau_2 ] .
\]
As before, the $M$-$N$ system can be obtained from the $M_+$-$N$
system by dividing out the cyclic symmetry carried by
$[\canr^\ext]$. In terms of graphs, the cyclic $\bbZ_3$ symmetry
corresponds to the three wings of the graph $\cE^{(12)}_1$ which
are transformed into each other by translation through the
$[\tau_j]$'s, and dividing out this symmetry gives exactly the
graph $\cE^{(12)}_2$ as the wings are identified whereas each
vertex on the middle axis splits into three nodes of identical
Perron-Frobenius weight. This way we understand the graph
$\cE^{(12)}_2$ as the label for the conjugation invariant
$Z_{\cE_2^{(12)}}$ of $Z_{\cE_1^{(12)}}$ which accidentally
happens to be the same as the self-conjugate $Z_{\cE^{(12)}}$
when specialized to $\SUd_9$ variables.
So here we have found some kind of two-fold degeneracy of the
modular invariant $Z_{\cE^{(12)}}$.

An even higher degeneracy appears for the modular invariant
of $\SUd_3$ which comes from the conformal embedding
$\SUd_3\subset{\mathit{SO}}(8)_1$.
Let $N\subset M_+$ be the corresponding local subfactor,
as usual. Let us briefly recall some facts about the
ambient ${\mathit{SO}}(8)_1$ theory.
It has four sectors, the basic ($0$),
vector (v), spinor (s) and conjugate spinor (c) module.
The conformal dimensions are given as $h_0=0$,
$h_\rmv=h_\rms=h_\rmc=1/2$, and the sectors obey
$\bbZ_2\times\bbZ_2$ fusion rules.
The Kac-Peterson matrices are given explicitly as
\[
S = \frac 12 \left( \begin{array}{rrrr}
1 & 1 & 1 & 1 \\ 1 & 1 & -1 & -1 \\
1 & -1 & 1 & -1 \\ 1 & -1 & -1 & 1
\end{array} \right) , \qquad
T = \E^{\pi\I/3} \left( \begin{array}{rrrr}
1 & 0 & 0 & 0 \\ 0 & -1 & 0 & 0 \\
0 & 0 & -1 & 0 \\ 0 & 0 & 0 & -1
\end{array} \right) .
\]
It is easy to see that there are exactly six
modular invariants, namely the six permutations
of v, s and c. Thus we are exclusively dealing with
automorphism invariants here, and in terms of
$\SUd_3$ variables they all specialize to the
same modular invariant
\begin{equation}
\label{ZcD6}
\cZ_{\mathcal{D}^{(6)}} =
\abs{\chi_{0,0} + \chi_{3,0} + \chi_{0,3}}^2
+ 3\abs{\chi_{1,1}}^2
\end{equation}
since the ${\mathit{SO}}(8)_1$ characters decompose upon
restriction to $\SUd$ variables into the level $3$
characters as
\[
\chi^0=\chi_{0,0} + \chi_{3,0} + \chi_{0,3},\qquad
\chi^\rmv=\chi^\rms=\chi^\rmc=\chi_{1,1}.
\]

The $\bbZ_2\times\bbZ_2$ fusion rules for
${\mathit{SO}}(8)_1$ models
were proven in the DHR framework in \cite{B2}, and
together with the conformal spin and statistics theorem
\cite{FG1,FRS2,GL2} we conclude that there is a system
$\{\id,\tau_\rmv,\tau_\rms,\tau_\rmc\}\subset\End(M_+)$
of braided endomorphisms, such that the statistics
S- and T-matrices are given exactly as above.
Because the statistics phases are given as
$\omega_\rmv=\omega_\rms=\omega_\rmc=-1$,
we can assume that the morphisms in the system
obey the $\bbZ_2\times\bbZ_2$ fusion rules even
by individual multiplication,
\[ \tau_\rmv^2=\tau_\rms^2=\tau_\rmc^2=\id\,, \qquad
\tau_\rmv \tau_\rms = \tau_\rms \tau_\rmv = \tau_\rmc \,,\]
thanks to Rehren's lemma \cite{R3}.
Hence we can extend $M_+$ in three ways as crossed products
by $\bbZ_2$, and the corresponding dual canonical endomorphism
sectors $[\canr]$ are respectively $[\id]\oplus[\tau_\rmv]$,
$[\id]\oplus[\tau_\rms]$ and $[\id]\oplus[\tau_\rmc]$.
We can also extend by the full group $\bbZ_2\times\bbZ_2$
giving instead
$[\id]\oplus[\tau_\rmv]\oplus[\tau_\rms]\oplus[\tau_\rmc]$.
Checking
$\lan\iota\la,\iota\mu\ran=\lan\canr\la,\mu\ran$
for $\la,\mu=\id,\tau_\rmv,\tau_\rms,\tau_\rmc$,
we find that there are only two $M$-$N$ sectors
for the $\bbZ_2$ extensions and only a single one,
namely $[\iota]$, for the full $\bbZ_2\times\bbZ_2$
extension. So $\tr Z=\#\MXN$ tells us that the
modular invariants of ${\mathit{SO}}(8)_1$ which
arise from the $M_+\subset M_+\rtimes\bbZ_2$ extensions
can only be the transpositions whereas the full
$M_+\subset M_+\rtimes(\bbZ_2\times\bbZ_2)$ extension
must produce one of the two cyclic permutations.
In fact it can produce both of them because the
matrices are relatively transpose and therefore
one can be obtained from the other by exchanging
braiding and opposite braiding.
It is easy to see that the extension by $\tau_\rmv$
gives exactly the transposition fixing v (and
analogously for s and c): By
$\langle\a^+_\la,\a^-_\mu\rangle\le\langle\canr\la,\mu\rangle$
we find $Z_{\rmv,\rms}=0=Z_{\rmv,\rmc}$ for
$[\canr]=[\id]\oplus[\tau_\rmv]$, leaving
only this possibility.

So what does the braided subfactor $N\subset M$
with $M$ one of these extensions give? Clearly,
they all produce the invariant of \erf{ZcD6}.
What are the relevant $\SUd$ graphs? A little
calculation shows easily that the $\bbZ_2$
extensions give the graph $\cD^{(6)}{}^*$ of
the Di Francesco-Zuber list whereas
$N\subset M_+\rtimes(\bbZ_2\times\bbZ_2)$
yields the same graph as $N\subset M_+$,
namely $\cD^{(6)}$. This is one of the examples
where different 2D CFT's, namely with different
extended coupling matrices (here the trivial one and
a non-trivial cyclic permutation of v, s, c),
are associated with the same graph.
Therefore we should not consider graphs as complete
labels of 2D extensions of some given chiral algebra.

By the way, a similar phenomenon seems already
to happen for the A-D-E invariants of $\SUz_k$.
For any D$_{\mathrm{even}}$ invariant, the exchange of
the two (``marked'') vertices sitting on the short legs
of the D-graph is a fusion rule automorphism of the
neutral system (see \cite[Subsect.\ 3.4]{BE2})
leaving conformal dimensions invariant,
but the exchange is invisible on the level of $\SUz_k$
characters since both correspond to the
character $\chi_{k/2}$.
So in this sense, all the D$_{\mathrm{even}}$
invariants are two-fold degenerate. A special case
is D$_{10}$ where one can also permute these two
vertices with the third vertex on the longest leg.
So here acts\footnote{This was pointed out to us
by K.-H.\ Rehren.} again the permutation group $S_3$.
But the four permutations involving this vertex are
visible on the level of $\SUz_{16}$ characters
(since $\chi_{k/2}$ and $\chi_2$ are different),
and they indeed produce the E$_7$ invariant which
is henceforth four-fold degenerate.
The A$_{k+1}$, E$_6$ and E$_8$ invariants do not
have a degeneracy.



\newcommand\bitem[2]{\bibitem{#1}{#2}}

\def\aam              {Acta Appl.\ Math. }
\def\aip              {Ann.\ Inst.\ H.\ Poincar\'e (Phys.\ Th\'eor.) }
\def\cmp              {Com\-mun.\ Math.\ Phys. }
\def\duke             {Duke Math.\ J. }
\def\ijm              {Intern.\ J. Math. }
\def\jfa              {J.\ Funct.\ Anal. }
\def\jmp              {J.\ Math.\ Phys. }
\def\lmp              {Lett.\ Math.\ Phys. }
\def\rmp              {Rev.\ Math.\ Phys. }
\def\inv              {Invent.\ Math. }
\def\mpl              {Mod.\ Phys.\ Lett. }
\def\nup              {Nucl.\ Phys. }
\def\nupp             {Nucl.\ Phys.\ (Proc.\ Suppl.) }
\def\adma             {Adv.\ Math. }
\def\physa            {Physica \textbf{A} }
\def\ijmp             {Int.\ J.\ Mod.\ Phys. }
\def\jp               {J.\ Phys. }
\def\fdp              {Fortschr.\ Phys. }
\def\pl               {Phys.\ Lett.}
\def\rims             {Publ.\ RIMS, Kyoto Univ. }



\end{document}